\documentclass[aop,seceqn,citesort,MSNbibl,dvips]{arximspdf}

\doi{10.1214/10-AOP633}
\volume{40}
\issue{2}
\pubyear{2012}
\firstpage{858}
\lastpage{889}

\makeatletter

\newcommand{\iint}{\int\!\!\int}

\newtheorem{theorem}{Theorem}[section]
\newtheorem{lemma}{Lemma}[section]
\newtheorem{prop}{Proposition}[section]

\newproclaim{remark}{Remark}[section]

\newtheorem{cor}{Corollary}[section]

\newcommand{\form}{\mathcal{E}}
\newcommand{\dom}{\mathcal{F}}
\newcommand{\rd}{{\mathbb R}^d}
\newcommand{\el}{\mathcal{L}}
\newcommand{\D}{\mathcal{D}}
\newcommand{\di}{\operatorname{diag}}

\makeatother

\begin{document}
\begin{frontmatter}

\title{Jump-type Hunt processes generated by lower bounded semi-Dirichlet forms}
\runtitle{Lower bounded semi-Dirichlet form}

\begin{aug}
\author[A]{\fnms{Masatoshi} \snm{Fukushima}\ead[label=e1]{fuku2@mx5.canvas.ne.jp}} and
\author[A]{\fnms{Toshihiro} \snm{Uemura}\corref{}\ead[label=e2]{t-uemura@kansai-u.ac.jp}}
\runauthor{M. Fukushima and T. Uemura}
\affiliation{Osaka University and Kansai University}
\address[A]{Branch of Mathematical Science\\
Faculty of Engineering Science\\
Osaka University\\
Toyonaka, Osaka 560-8531\\
Japan \\
and\\
Department of Mathematics\\
Faculty of Engineering Science\\
Kansai University\\
Suita, Osaka 564-8680\\
Japan \\
\printead{e1}\\
\phantom{E-mail: }\printead*{e2}} 
\end{aug}

\received{\smonth{3} \syear{2010}}
\revised{\smonth{11} \syear{2010}}

%
\begin{abstract}
Let $E$ be a locally compact separable metric space and $m$ be a~%
positive Radon measure on it. Given a nonnegative function $k$ defined on
$E\times E$ off the diagonal whose anti-symmetric part is assumed to be less
singular than the symmetric part, we construct an associated regular
lower bounded
semi-Dirichlet form $\eta$ on $L^2(E;m)$ producing a Hunt process~%
$X^0$ on $E$
whose jump behaviours are governed\vspace*{1pt} by $k$. For an arbitrary open subset
$D\subset E$, we also construct a Hunt process $X^{D,0}$ on $D$ in an
analogous manner.
When $D$ is relatively compact, we show that $X^{D,0}$ is censored in
the sense
that it admits no killing inside $D$ and killed only when the path
approaches to
the boundary. When $E$ is a $d$-dimensional Euclidean space and $m$ is
the Lebesgue
measure, a typical example of $X^0$ is the stable-like process that
will be also
identified with the solution of a martingale problem up to an $\eta
$-polar set of starting
points. Approachability\vspace*{1pt} to the boundary $\partial D$ in finite time of
its censored
process $X^{D,0}$ on a~bounded open subset~$D$ will be examined in
terms of the
polarity of $\partial D$ for the symmetric stable processes with
indices that bound
the variable exponent $\alpha(x)$.
\end{abstract}

%

\end{frontmatter}

\section{Introduction}\label{sec1}
Let $E$ be a locally compact separable metric space equipped with a
metric $d$,
$m$ be a positive Radon measure with full topological support and
$k(x,y)$ be
a nonnegative Borel measurable function on the space $E\times
E\setminus\di$,
where $\di$ denotes the diagonal set $\{(x,x) \dvtx x \in E\}$. A purpose
of the present
paper is to construct Hunt processes on $E$ and on its subsets with
jump behaviors
being governed by the kernel $k$ by using general results on a lower bounded
semi-Dirichlet form on $L^2(E;m)$.

The inner product and the norm in $L^2(E;m)$ are denoted by $(\cdot
,\cdot)$
and \mbox{$\Vert\cdot\Vert$}, respectively. Let $\dom$ be a dense linear
subspace of
$L^2(E;m)$ such that $u\wedge1\in\dom$ whenever $u\in\dom$. A (not
necessarily
symmetric) bilinear form $\eta$ on~$\dom$ is called\vadjust{\goodbreak} a \textit{lower
bounded closed
form} if the following three conditions are satisfied: we set $
{\eta_\beta(u,v)
=\eta(u,v)+\beta(u,v), u,v\in\dom}$.
There exists a \mbox{$\beta_0\ge0$} such that:

\begin{longlist}[(B.2)]
\item[(B.1)] (lower boundedness); for any $u\in\dom$,
$\eta_{\beta_0}(u,u) \ge0$.

\item[(B.2)] (sector condition); for any $u,v\in\dom$,
\[
| \eta(u,v)| \le K
\sqrt{\eta_{\beta_0}(u,u)} \cdot\sqrt{\eta_{\beta_0}(v,v)}
\]
for some constant $K\ge1$.

\item[(B.3)] (completeness); the space $\dom$ is complete with
respect to the norm
$\eta^{1/2}_{\alpha}(\cdot, \cdot)$ for some, or equivalently, for all
$\alpha>\beta_0$.
\end{longlist}

For a lower bounded closed form $(\eta,\dom)$ on $L^2(E;m)$, there
exist unique
semigroups $\{T_t;t>0\}, \{\widehat T_t;t>0\}$ of linear operators on
$L^2(E;m)$ satisfying
%
%
\begin{eqnarray}\label{equ1.1}
(T_tf,g)=(f,\widehat T_t g),\nonumber\\[-8pt]\\[-8pt]
&&\eqntext{f,g \in L^2(E;m), \Vert T_t\Vert\le
e^{\beta_0t},
\Vert\widehat T_t\Vert\le e^{\beta_0t}, t>0,}
\end{eqnarray}
such that their Laplace transforms $G_\alpha$ and $\widehat G_\alpha$ are
determined for $\alpha>\beta_0$ by
\begin{eqnarray*}
&G_\alpha f, \widehat G_\alpha f\in\dom,\qquad \eta_\alpha(G_\alpha f,u)
=\eta_\alpha(u, \widehat G_\alpha f)=(f,u),&\\
&f\in L^2(E;m),
u\in\dom.&
\end{eqnarray*}
See the first part of Section \ref{sec3} for more details. $\{T_t;t>0\}
$ is said
to be
\textit{Markovian} if $0\le T_t f \le1, t>0$, whenever $f\in L^2(E;m),
0\le
f\le1$. It was shown by Kunita \cite{K69} that the semigroup $\{
T_t;t>0\}$
is Markovian if and only if
%
%
\begin{equation}\label{equ1.2}
Uu\in\dom\quad\mbox{and}\quad \eta(Uu,u-Uu) \ge0 \qquad\mbox{for any }
u\in\dom,
\end{equation}
where $Uu$ denotes the unit contraction of $u$: $Uu=(0\vee u)\wedge1$.
A lower
bounded closed form $(\eta, \dom)$ on $L^2(E;m)$ satisfying (\ref{equ1.2})
will be called a~\textit{lower bounded semi-Dirichlet form} on $L^2(E;m)$. The term
``semi'' is
added to indicate that the dual semigroup $\{\widehat T_t; t>0\}$ may
not be
Markovian although it is positivity preserving. As we shall see in
Section \ref{sec3} for a
lower bounded semi-Dirichlet form $\eta$ which is regular in the sense
stated below,
if the associated dual semigroup $\{\widehat T_t; t>0\}$ were
Markovian, or equivalently,
if $m$ were excessive, then $\eta$ is necessarily a nonnegative
definite closed
form, namely,~$\beta_0$ in conditions (B.1), (B.3) [resp.,
(B.2)] can
be retaken to be $0$ (resp., $1$).

A lower bounded semi-Dirichlet form $(\eta, \dom)$ is said to be \textit
{regular} if
$\dom\cap C_0(E)$ is uniformly dense in $C_0(E)$ and $\eta_\alpha
$-dense in $\dom$
for $\alpha>\beta_0$, where~$C_0(E)$ denotes the space of continuous
functions on
$E$ with compact support. Carrillo-Menendez \cite{C75} constructed
a Hunt process
properly associated with any regular lower bounded semi-Dirichlet form on
$L^2(E;m)$ by reducing the situation to the case where $\eta$ is nonnegative
definite. We shall show in Section \ref{sec4} that a direct
construction is
possible without such a reduction.\vadjust{\goodbreak}

Later on, the nonnegative definite semi-Dirichlet form was investigated by
Ma, Oberbeck and R\"{o}ckner \cite{MOR94} and
Fitzsimmons \cite{Fi01}
specifically in a~general context of the quasi-regular Dirichlet form
and the
special standard process. However, in producing the forms $\eta$ from
nonsymmetric kernels $k$ corresponding to a considerably wide class of jump
type Hunt processes in finite dimensions whose dual semigroups need not be
Markovian, we will be forced to allow positive $\beta_0$.

To be more precise, we set for $x,y \in E, x\not=y$,
%
%
\begin{equation}\label{kernel}
k_s(x,y):=\tfrac12 \{ k(x,y)+k(y,x)\} \quad\mbox{and}\quad
k_a(x,y):=\tfrac12 \{ k(x,y)-k(y,x)\},\hspace*{-36pt}
\end{equation}
that is, the kernel $k_s(x,y)$ denotes the symmetrized one of $k$, while
$k_a(x,y)$~re\-presents the anti-symmetric part of $k$. We impose four conditions
\mbox{(\ref{cond-1})--(\ref{cond-4})} on $k_s$ and $k_a$ stated below. Condition
(\ref{cond-1}) on $k_s$ is nearly optimal for us to work with the symmetric
Dirichlet form (\ref{symmetric-form}) defined below, while
condi\-tions~(\ref{cond-2})--(\ref{cond-4}) require $k_a$ to be less singular than $k_s$.

Let   conditions (\ref{cond-1})--(\ref{cond-4}) be in force on $k$.
Denote by $C_0^{\mathrm{lip}}(E)$
the space of uniformly Lipschitz continuous functions
on $E$ with compact support. We also let
%
%
\begin{equation}\label{symmetric-form}
\cases{
\displaystyle \form(u,v):= \iint_{E\times E\setminus\di}
\bigl(u(y)-u(x)\bigr)
\bigl(v(y)-v(x)\bigr)\vspace*{2pt}\cr
\hspace*{98pt}{}\times k_s(x,y)m(dx)m(dy), \vspace*{2pt}\cr
\displaystyle \dom^r=\{u\in L^2(E;m)\dvtx u \mbox{ is Borel measurable and }
\form(u,u)<\infty\}.}\hspace*{-35pt}
\end{equation}
$(\form, \dom^r)$ is a symmetric Dirichlet form on $L^2(E;m)$ and
$\dom^r$ contains
the space~$C_0^{\mathrm{lip}}(E)$. We denote by $\dom^0$ the $\form
_1$-closure of
$C_0^{\mathrm{lip}}(E)$ in $\dom^r$. $(\form, \dom^0)$ is then a~regular Dirichlet
form on $L^2(E;m)$ (cf. \cite{FOT94}, Example 1.2.4, Theorem~3.1.1 and see also \cite{U02} and \cite{U04}).

For $u\in C_0^{\mathrm{lip}}(E)$ and $n\in{\mathbb N}$, the integral
%
%
\begin{equation}\label{operator}\quad
\mathcal{L}^nu(x):=\int_{\{y\in E\dvtx d(x,y)>1/n\}} \bigl(u(y)-u(x)
\bigr)k(x,y)m(dy),\qquad
x\in E,
\end{equation}
makes sense. We prove in Proposition \ref{prop2.1} and Theorem \ref
{theo2.1} in Section \ref{sec2}
that the finite limit
%
%
\begin{equation}\label{semi-form}
\eta(u,v)=-\lim_{n\to\infty}\int_E \mathcal{L}^nu(x) v(x)
m(dx)\qquad \mbox{for }
u,v\in C_0^{\mathrm{lip}}(E),
\end{equation}
exists, $\eta$ extends to $\dom^0\times\dom^0$ and $(\eta,\dom
^0)$ is a lower
bounded semi-Dirichlet form on $L^2(E;m)$ with parameter $\beta
_0=8(C_1\vee C_2C_3)
(\ge0)$ where $C_1$--$C_3$ are constants appearing in   conditions
(\ref{cond-2})--(\ref{cond-4}). Furthermore, the form $\form$ is
shown to be a
\textit{reference} (\textit{symmetric Dirichlet}) \textit{form} of $\eta
$ in
the sense
that, for each fixed $\alpha>\beta_0$,
%
%
\begin{equation}\label{ref-form}
c_1 \form_1(u,u)\le\eta_\alpha(u,u)\le c_2 \form_1(u,u),\qquad
u\in\dom^0,
\end{equation}
for some positive constants $c_1, c_2$ independent of $u\in\dom^0$.
Therefore,
$(\eta, \dom^0)$ becomes a regular\vadjust{\goodbreak} lower bounded semi-Dirichlet form
on $L^2(E;m)$
and gives rise to an associated Hunt process $X^0=(X_t^0, P_x^0)$ on $E$.
We call $X^0$ the \textit{minimal Hunt process} associated with the form
$\eta$.
Equation (\ref{semi-form}) indicates that the limit of $\mathcal{L}^n$
in $n$
plays a role of
a pre-generator of $X^0$ informally.

If we define the kernel $k^*$ by
%
%
\begin{equation}\label{reverse}
k^*(x,y):=k(y,x),\qquad x,y\in E, x\not=y,
\end{equation}
and the form $\eta^*$ by (\ref{operator}) and (\ref{semi-form}) with
$k^*$ in
place of $k$, we have the same conclusions as above for $\eta^*$
(Corollary \ref{cor2.1}
of Section \ref{sec2}). In particular, there exists a minimal Hunt process
$X^{0*}$ associated
with the form $\eta^*$.

In the second half of Section \ref{sec3}, we are concerned with a
killed dual semigroup
$\{ e^{-\beta t} \widehat{T_t}; t>0\}$, which can be verified to be
Markovian for a large $\beta>0$ but only for a restricted subfamily of the
forms $\eta$ considered in Section \ref{sec2} (lower order cases). For
a higher
order $\eta$,
the killed dual semigroup may not be Markovian no matter how big $\beta
$ is.
We shall also exhibit an example of a one-dimensional probability kernel
$k$ [$\int_{{\mathbb R}^1}k(x,y) \,dy=1$] with $m$ being the Lebesgue measure,
for which the associated semi-Dirichlet form $\eta$ is not nonnegative definite
and accordingly the associated dual semigroup itself is non-Markovian.

When $E={\mathbb R}^d$ the $d$-dimensional Euclidean space and
$m(dx)=dx$ the~Le\-besgue measure on it, we shall verify in Section \ref{sec5} that our
requirements~(\ref{cond-1})--(\ref{cond-4}) on the kernel $k(x,y)$ are fulfilled by
%
%
\begin{eqnarray}\label{stable-like-kernel}
k^{(1)}(x,y)&=&w(x)|x-y|^{-d-\alpha(x)},\nonumber\\[-8pt]\\[-8pt]
k^{(1)*}(x,y)&=&w(y)|x-y|^{-d-\alpha(y)},\qquad
x,y \in{\mathbb R}^d, x\neq y,\nonumber
\end{eqnarray}
for $w(x)$ given by (\ref{weight}) and $\alpha(x)$ satisfying the
bounds (\ref{stable-like}).
A Markov process corresponding to $k^{(1)}$ is called a \textit
{stable-like process}
and has been constructed by Bass \cite{B88a} as a unique solution
to a martingale
problem. In this case, we shall prove that the minimal Hunt process
associated with
the corresponding form $\eta$ is conservative and actually a solution
to the same
martingale problem, identifying it with the one constructed in \cite
{B88a} up to an
$\eta$-polar set of starting points.

In Section \ref{sec6}, we consider an arbitrary open subset $D$ of $E$. Define
$m_D$ by
$m_D(B)=m(B\cap D)$ for any Borel set $B\subset E$. By replacing $E$
and $m$
with~$D$ and $m_D$, respectively, in (\ref{symmetric-form}), we obtain
a symmetric
Dirichlet form $(\form_D, \dom_D^r)$ on $L^2(D;m_D)$. Denote by
$\overline D$ the
closure of $D$ and by $C_0^{\mathrm{lip}}(\overline D)$ the restriction to
$\overline D$
of the space $C_0^{\mathrm{lip}}(E)$. We also denote by $C_0^{\mathrm{lip}}(D)$
the space of
uniformly Lipschitz continuous functions on $D$ with compact support in $D$.
Let $\dom_{\bar D}$ and $\dom_D^0$ be the $\form_{D,1}$-closures of
$C_0^{\mathrm{lip}}(\overline D)$ and $C_0^{\mathrm{lip}}(D)$,
respectively, in
$\dom_D^r$.
Then $(\form_D, \dom_{\bar D})$ is a regular symmetric Dirichlet form on
$L^2(\overline D;m_D)$, while $(\form_D^0, \dom_D^0)$ is a regular symmetric
Dirichlet form on $L^2(D;m_D)$ where~$\form_D^0$ is the restriction of
$\form_D$
to $\dom_D^0\times\dom_D^0$.

By making the same replacement in (\ref{operator}) and (\ref
{semi-form}), we get a
form $\eta_D$ on $C_0^{\mathrm{lip}}(\overline D)\times C_0^{\mathrm
{lip}}(\overline D)$,
which extends\vadjust{\goodbreak} to $\dom_{\bar D}\times\dom_{\bar D}$ to be a regular
lower bounded
semi-Dirichlet form on $L^2(\overline D;m_D)$ possessing $\form_D$ as
its reference
form, yielding an associated Hunt process $X^{\bar D}$ on $\overline D$.
We also consider the restriction $\eta_D^0$ of $\eta_D$ to $\dom
_D^0\times\dom_D^0$
so that $(\eta_D^0, \dom_D^0)$ is a regular lower bounded
semi-Dirichlet form on
$L^2(D;m_D)$ possessing $\form_D^0$ as its reference form. We shall
show in Section \ref{sec6}
that the part process $X^{D,0}$ of $X^{\bar D}$ on $D$, namely, the
Hunt process
obtained from $X^{\bar D}$ by killing upon hitting the boundary
$\partial D$,
is properly associated with $(\eta_D^0, \dom_D^0)$.

We shall also prove in Section \ref{sec6} that $X^{\bar D}$ admits no
jump from $D$ to~$\partial D$, and furthermore when $D$ is relatively compact,
$X^{\bar D}$ is conservative so that $X^{D,0}$ admits no killing inside $D$ and its
sample path
is killed only when it approaches to the boundary $\partial D$.
$X^{D,0}$ is
accordingly different from the part process of $X^0$ on the set $D$ in general
because the sample path of $X^0$ may jump from $D$ to $E\setminus D$
resulting in a killing inside $D$ of its part process. By adopting
$k^*$ instead
of $k$, we get in an analogous manner Hunt processes~$X^{\bar D*}$ on
$\overline D$
and $X^{D,0*}$ on $D$ satisfying the same properties as above.

When $(\form, \dom^r)$ is the Dirichlet form on $L^2(\rd)$ of a
symmetric stable
process on $\rd$, the space $\dom^0$ is identical with $\dom^r$. In
this case,
for an arbitrary open set $D\subset\rd$, the symmetric Hunt process
on $D$
associated with $(\form_D^0, \dom_D^0)$ is a \textit{censored stable
process} on $D$
in the sense of Bogdan, Burdzy and Chen~\cite{BBC03}. It
was further
shown in \cite{BBC03} that, if $D$ is a $d$-set, then the space~$\dom
_{\bar D}$
coincides with $\dom^r_D$ so that the symmetric Hunt process on
$\overline D$
associated with $(\form_D, \dom^r_D)$ was called a \textit{reflecting
stable process}
over $\overline D$.

For the nonsymmetric kernel $k^{(1)}$ on $\rd$ as (\ref{stable-like-kernel}),
associated Hunt processes~$X^{D,0}, X^{D,0*}$ on an arbitrary open set
$D\subset\rd$ may well be called \textit{censored stable-like processes}
in view of
the stated properties of them. However, it is harder in this case to
identify the
space $\dom_{\bar D}$ with $\dom^r_D$, and accordingly we call the
associated Hunt
processes $X^{\bar D}, X^{\bar D*}$ over $\overline D$ \textit{modified
reflecting
stable-like processes} analogously to the Brownian motion case (cf.
\cite{F99}).
At the end of Section \ref{sec6}, we give sufficient conditions in
terms of the
upper and lower
bounds of the variable exponent $\alpha(x)$ for the approachability in
finite time
of the censored stable-like processes to the boundary.

We are grateful to Professor Yoichi Oshima for providing us with his unpublished
lecture notes \cite{O88} on nonsymmetric Dirichlet forms as well as
an updated
version of a part of it, which are very valuable for us.\vspace*{-3pt}

\section{Construction of a lower bounded semi-Dirichlet form from
$k$}\label{sec2}

Throughout this section, we make the following assumptions on a nonnegative
Borel measurable function $k(x,y)$ on $E\times E\setminus\di$:
%
%
\begin{eqnarray}\quad
\label{cond-1}
&\displaystyle M_s\in L^2_{\mathrm{loc}}(E;m) \qquad\mbox{for } M_s(x)= \int_{y\not=x}
\bigl( 1\wedge d(x,y)^2 \bigr) k_s(x,y)m(dy),&\nonumber\\[-8pt]\\[-8pt]
&&\eqntext{x\in E,}
\\
%
\label{cond-2}
&\displaystyle C_1:=\sup_{x\in E} \int_{d(x,y)\ge1} |k_a(x,y)|m(dy) <\infty,&
\end{eqnarray}
and there exists a constant $\gamma\in(0,1]$ such that
%
%
\begin{equation}\label{cond-3}
C_2:= \sup_{x\in E} \int_{d(x,y)<1}
| k_a(x,y)|^{\gamma} m(dy)<\infty,
\end{equation}
and furthermore,
for some constant $C_3\ge0$,
%
%
\begin{eqnarray}\label{cond-4}
| k_a(x,y)|^{2-\gamma} \le C_3 k_s(x,y)\qquad \mbox
{for any } x,y\in E \nonumber\\[-8pt]\\[-8pt]
&&\eqntext{\mbox{with } 0<d(x,y) \le1.}
\end{eqnarray}

For each $n\in\mathbb N$, define ${\mathcal L}^nu$ for $u\in
C_0^{\mathrm{lip}}(E)$ by
(\ref{operator}) and $\eta^n(u,v)$ for $u,v \in C_0^{\mathrm{lip}}(E)$ by
%
%
\begin{equation}\label{form}
\eta^n(u,v):=-\int_E \mathcal{L}^nu(x) v(x) m(dx),
\end{equation}
the integral on the right-hand side being absolutely convergent by
(\ref{cond-1}).~We note that any $u\in C_0^{\mathrm{lip}}(E)$ belongs to the domain
$\dom
^r$ of the form
$\form$ defined by~(\ref{symmetric-form}). In fact, if we denote by
$K$ the support
of $u$, then $\form(u,u)$ is dominated by twice the integral of
$(u(x)-u(y))^2k_s(x,y)m(dx)m(dy)$ on $K\times E$, which is finite by
(\ref{cond-1}).

$\form(u,v)$ admits also an alternative expression for $u,v \in
C_0^{\mathrm{lip}}(E)$,
\[
\form(u,v):=\iint_{E\times E\setminus\di} \bigl(u(y)-u(x)\bigr)
\bigl(v(y)-v(x)\bigr) k(x,y) m(dx)m(dy),
\]
because the right-hand side of the above can be seen to be equal to the same
integral with $k(y,x)$ in place of $k(x,y)$ by interchanging the
variables~$x, y$, and we arrive at the expression in (\ref{symmetric-form}) by
averaging.
In particular, $ \form(u,v)=\lim_{n\to\infty}\form
^n(u,v)$ for
$u,v\in C_0^{\mathrm{lip}}(E)$ where
%
%
\begin{equation}\label{form-1}
\form^n(u,v):= \iint_{d(x,y)> 1/n} \bigl(u(y)-u(x)\bigr)
\bigl(v(y)-v(x)\bigr)k(x,y)m(dx)m(dy).\hspace*{-35pt}
\end{equation}
%

\begin{prop}\label{prop2.1} Assume (\ref{cond-1})--(\ref{cond-4}). Then
for all
$u,v\in C_0^{\mathrm{lip}}(E)$, the limit
\[
\eta(u,v)=\lim_{n\to\infty} \eta^n(u,v)
\]
exists. Moreover, the limit has the following expression:
%
%
\begin{equation}\label{3.4}\qquad
\eta(u,v)= \frac12 \form(u,v) + \iint_{y\not=x}
\bigl(u(x)-u(y)\bigr) v(y)
k_a(x,y)m(dx)m(dy),
\end{equation}
where $\form$ is defined by (\ref{symmetric-form}) and the integral
on the right-hand side is absolutely convergent.
\end{prop}
\begin{pf}
For $u,v \in C_0^{\mathrm{lip}}(E)$, we have
\begin{eqnarray*}
\eta^n(u,v)-\eta^n(v,u) &=&
-\iint_{d(x,y)> 1/n} \bigl(u(y)-u(x)\bigr)v(x) k(x,y)m(dx)m(dy) \\
&&{} + \iint_{d(x,y)> 1/n}\bigl(v(y)-v(x)\bigr) u(x)
k(x,y)m(dx)m(dy) \\
&=& -\iint_{d(x,y)> 1/n}u(y)v(x) k(x,y)m(dx)m(dy) \\
&&{} + \iint_{d(x,y)> 1/n}v(y)u(x) k(x,y)m(dx)m(dy) \\
&=& 2 \iint_{d(x,y)> 1/n} u(x) v(y) k_a(x,y)m(dx)m(dy),
\end{eqnarray*}
and further
\begin{eqnarray*}
&&\eta^n(u,v) +\eta^n(v,u)\\[-0.75pt]
&&\qquad= - \iint_{d(x,y)\ge1/n}\bigl(u(y)-u(x)\bigr) v(x)
k(x,y)m(dx)m(dy) \\[-0.75pt]
&&\qquad\quad{} -\iint_{d(x,y)\ge1/n}\bigl(v(y)-v(x)\bigr) u(x)
k(x,y)m(dx)m(dy) \\[-0.75pt]
&&\qquad= \iint_{d(x,y)\ge1/n}\bigl(u(y)-u(x)\bigr)
\bigl(v(y)-v(x)\bigr)
k(x,y)m(dx)m(dy) \\[-0.75pt]
&&\qquad\quad{} -\iint_{d(x,y)\ge1/n}\bigl(u(y)-u(x)\bigr) v(y)
k(x,y)m(dx)m(dy) \\[-0.75pt]
&&\qquad\quad{} -\iint_{d(x,y)\ge1/n} \bigl(v(y)-v(x)
\bigr)u(x)k(x,y)m(dx)m(dy) \\[-0.75pt]
&&\qquad= \form^n(u,v) -2\iint_{d(x,y)\ge1/n}
u(y)v(y)k_a(x,y)m(dx)m(dy).
\end{eqnarray*}
By adding up the obtained identities, we get for $u,v\in C_0^{\mathrm{lip}}(E)$,
%
%
\begin{eqnarray}
\label{(8)}
2 \eta^n(u,v) &=& \form^n(u,v) + 2 \iint_{d(x,y)> 1/n}
\bigl(u(x)-u(y)\bigr) v(y)\nonumber\\[-8pt]\\[-8pt]
&&\hspace*{112pt}{} \times
k_a(x,y) m(dx)m(dy).\nonumber
\end{eqnarray}

Since $\form^n(u,v)$ converges to $\form(u,v)$ as $n\to\infty$,
it remains to see that the second term of the right-hand side also
converges absolutely as $n\to\infty$ for each $u,v \in C_0^{\mathrm{lip}}(E)$.

From the Schwarz inequality and (\ref{cond-2}),
we see that
\begin{eqnarray*}
&&\iint_{d(x,y)> 1/n} \bigl|
\bigl(u(x)-u(y)\bigr) v(y)
k_a(x,y) \bigr| m(dx)m(dy) \\[1pt]
&&\qquad\le\iint_{1/n<d(x,y)< 1} |u(x)-u(y)| \cdot
|v(y)|
|k_a(x,y)|^{\gamma/2}\\[1pt]
&&\qquad\quad\hspace*{67pt}{}\times|k_a(x,y)|^{1-\gamma/2}
m(dx)m(dy) \\[1pt]
&&\qquad\quad{} + \iint_{d(x,y)\ge1} |u(x)-u(y)| \cdot
|v(y)|
k_s(x,y)^{1/2}|k_a(x,y)|^{1/2}m(dx)m(dy) \\[1pt]
&&\qquad\le\sqrt{ \iint_{1/n<d(x,y)<1} \bigl(u(x)-u(y)\bigr)^2
|k_a(x,y)|^{2-\gamma} m(dx)m(dy) } \\[1pt]
&&\qquad\quad{}\times
\sqrt{\iint_{1/n<d(x,y)<1} v(y)^2 |k_a(x,y)|^{\gamma}
m(dx)m(dy)} \\[1pt]
&&\qquad\quad{} + \sqrt{C_1} \Vert v\Vert \sqrt{ \iint_{d(x,y)\ge1}
\bigl(u(x)-u(y)\bigr)^2 k_s(x,y)m(dx)m(dy)}.
\end{eqnarray*}
So, by making use of assumptions (\ref{cond-3}) and (\ref{cond-4}) and
an elementary inequality
$\sqrt{A}+\sqrt{B} \le\sqrt{2}\sqrt{A+B}$
holding for $A\ge0 $ and $B\ge0$,
we have
\begin{eqnarray*}
&&
\iint_{d(x,y)> 1/n}\bigl| \bigl(u(x)-u(y)\bigr) v(y) k_a(x,y) \bigr|
m(dx)m(dy) \\[2pt]
&&\qquad\le\sqrt{2} \sqrt{C_1\vee C_2 C_3}
\Vert v\Vert \cdot\sqrt{\form^n(u,u)}.
\end{eqnarray*}
Then taking $n\to\infty$,
\begin{eqnarray*}
&&\iint_{y\not=x} \bigl| \bigl(u(x)-u(y)\bigr) v(y) k_a(x,y) \bigr|
m(dx)m(dy) \\[2pt]
&&\qquad\le\sqrt{2} \sqrt{C_1\vee C_2 C_3}
\Vert v\Vert \cdot\sqrt{\form(u,u)}<\infty
\end{eqnarray*}
as was to be proved.
\end{pf}

For $u,v\in C_0^{\mathrm{lip}}(E)$, set
\[
\eta_{\beta}(u,v)=\eta(u,v)+\beta(u,v),\qquad \beta>0,
\]
and
%
%
\begin{equation}\label{B}
B(u,v):=\iint_{x\not=y} \bigl(u(x)-u(y)\bigr) v(y) k_a(x,y)m(dx)m(dy).
\end{equation}
Then   equation (\ref{3.4}) reads
%
%
\begin{equation}\label{(9)}
\eta(u,v)= \tfrac12 \form(u,v)+ B(u,v),\qquad u,v\in C_0^{\mathrm{lip}}(E),
\end{equation}
while we get from the proof of the preceding proposition
%
%
\begin{equation}\label{eqn:bound-B}
| B(u,v)| \le C_4 \Vert v\Vert \sqrt{\form(u,u)},
\end{equation}
where $C_4=\sqrt{2}\cdot\sqrt{C_1\vee C_2C_3}$. Now we put
$\beta_0:=4(C_4)^2=8 (C_1\vee C_2C_3)$.

From   equation (\ref{(9)}) and the bound (\ref{eqn:bound-B}),
we have for $ u\in C_0^{\mathrm{lip}}(E)$,
\begin{eqnarray*}
\eta_{\beta_0}(u,u)&=&
\tfrac14 \form_{\beta_0}(u,u)+\tfrac14 \form(u,u)+\tfrac34\beta
_0\Vert u\Vert^2+B(u,u)\\
&\ge&
\tfrac14 \form_{\beta_0}(u,u)+\sqrt{3} C_4 \sqrt{\form(u,u)} \Vert
u\Vert+B(u,u)
\ge\tfrac14\form_{\beta_0}(u,u).
\end{eqnarray*}

Further, for $u,v\in C_0^{\mathrm{lip}}(E)$,
\begin{eqnarray*}
| \eta(u,v)|
&\le& \tfrac12 | \form(u,v) | + |B(u,v)|\\
&\le& \tfrac12 \sqrt{\form(u,u)} \sqrt{\form(v,v)}
+ C_4 \Vert v\Vert \sqrt{\form(u,u)}
\\
&\le& \tfrac12 \bigl( \sqrt{\form(v,v)} + 2C_4 \Vert v\Vert
\bigr)
\sqrt{\form(u,u)} \\
&\le&
\tfrac{\sqrt{2}}{2} \sqrt{\form_{\beta_0}(v,v)}
\sqrt{\form_{\beta_0}(u,u)}.
\end{eqnarray*}

So it also follows that
%
%
\begin{equation}\label{sector}
| \eta(u,v) | \le2\sqrt{2} \sqrt{\eta_{\beta_0}(u,u)}
\sqrt{\eta_{\beta_0}(v,v)}
\end{equation}
and
%
%
\begin{equation}\label{eqn:comparison}
\tfrac{1}{4} \form_{\beta_0}(u,u) \le\eta_{\beta_0}(u,u) \le
\tfrac{2+\sqrt{2}}{2} \form_{\beta_0}(u,u),\qquad u,v \in C_0^{\mathrm{lip}}(E).
\end{equation}

Let $\dom^0$ be the $\form_1$-closure of $C_0^{\mathrm{lip}}(E)$ in
$\dom^r$.
Since $\dom^0$ is complete with respect to $\form_\alpha$ for any
$\alpha>0$, the estimates obtained above readily lead us to the first
conclusion of
the following theorem.
%
%
\begin{theorem}\label{theo2.1} Assume (\ref{cond-1})--(\ref{cond-4}).
Then the form $\eta$ defined by Proposition~\ref{prop2.1} extends from
$C_0^{\mathrm {lip}}(E)\times C_0^{\mathrm{lip}}(E)$ to
$\dom^0\times\dom^0$ to be a lower bounded closed form on $L^2(E;m)$
satisfying \textup{(B.1)--(B.3)} with $\beta_0=8(C_1\vee C_2C_3), K=2\sqrt{2}$
and possessing $(\form,\dom^0)$ as a reference form in the sense of
(\ref{ref-form}).

Furthermore, the pair $(\eta, \dom^0)$ is a regular lower bounded
semi-Dirichlet
form on $L^2(E;m)$.
\end{theorem}

We note that the above constant $\beta_0$ is equal to $0$ if $k$ is symmetric:
$k(x,y)=k(y,x), (x,y)\in E\times E\setminus\di$.
\begin{pf*}{Proof of Theorem \ref{theo2.1}}
It suffices to prove the contraction proper\-ty~(\ref{equ1.2}) for
the present
pair $(\eta, \dom^0)$. We first show this for $u \in
C_0^{\mathrm{lip}}(E)$.\vadjust{\goodbreak}
Note that $Uu \in C_0^{\mathrm{lip}}(E)$ and, for $n\in{\mathbb N}$,
\begin{eqnarray*}
&&\eta^n(Uu,u-Uu) \\
&&\qquad= -\iint_{d(x,y)>1/n} \bigl(Uu(y)-Uu(x)\bigr)
\bigl(u(x)-Uu(x)\bigr) k(x,y)m(dx)m(dy) \\
&&\qquad= \iint_{\{d(x,y)>1/n\}\cap\{x \dvtx u(x)\ge1\}}
\bigl(1- Uu(y)\bigr) \bigl( u(x)-1 \bigr)k(x,y)m(dx)m(dy) \\
&&\qquad\quad{} -\iint_{\{d(x,y)>1/n\}\cap\{x \dvtx u(x)\le0\}} Uu(y)
u(x) k(x,y)m(dx)m(dy) \\
&&\qquad\ge0.
\end{eqnarray*}
Then, we have by Proposition \ref{prop2.1}
\[
\eta(Uu,u-Uu) =\lim_{n\to\infty} \eta^n(Uu,u-Uu) \ge0.
\]

Following a method in \cite{MR92}, Lemma 4.9, we next prove (\ref
{equ1.2}) for any
$u\in\dom^0$. Choose a sequence $\{u_{\ell}\} \subset C_0^{\mathrm{lip}}(E)$
which is $\form_1$-convergent to $u$. Then
%
%
\begin{equation}\label{U1}
\Vert Uu_{\ell}-Uu\Vert\rightarrow0,\qquad \ell\to\infty,
\end{equation}
because $U$ is easily seen to be a continuous operator from $L^2(E;m)$ to
$L^2(E;m)$. Fix $\alpha>\beta_0$. We then get from (\ref{ref-form})
the boundedness
\[
\sup_\ell\eta_\alpha(Uu_\ell,Uu_\ell) \le C_2 \sup_\ell\form
_1(u_\ell,u_\ell)
<\infty.
\]
On the other hand, using the dual resolvent $\widehat G_\alpha$ associated
with the lower bounded closed form $(\eta,\dom^0)$, we see from
equation (\ref{equ3.1})
below that, for any $g\in L^2(E;m)$,
\[
\eta_\alpha(Uu_\ell, \widehat G_\alpha g)=(Uu_\ell, g)
\rightarrow
(Uu,g)=\eta_\alpha(Uu, \widehat G_\alpha g), \qquad\ell\to\infty.
\]
Since $\{\widehat G_\alpha g\dvtx g\in L^2(E,m)\}$ is $\eta_\alpha$-dense
in $\dom^0$,
we can conclude by making use of the above $\eta_\alpha$-bound of $\{
Uu_\ell\}$
and the sector condition (B.2) that $\{Uu_\ell\}$ is $\eta
_\alpha$-weakly
convergent to $Uu$ as $\ell\to\infty$. In particular, by the above
$\eta_\alpha$-bound and (B.2) again, we have
%
%
\begin{equation}\label{U2}
\eta_\alpha(Uu_\ell, u_\ell) \rightarrow\eta_\alpha(Uu,u),\qquad
\ell\to\infty.
\end{equation}

We consider the dual form $\widehat\eta$ and the symmetrizing form
$\tilde\eta$
of $\eta$ defined~by
\[
\widehat\eta(u,v)=\eta(v,u),\qquad \tilde\eta(u,v)=\tfrac12
\bigl(\eta(u,v)
+\eta(v,u)\bigr),\qquad u,v\in\dom^0.
\]
In the same way as above, we can see that $\{Uu_\ell\}$ converges as
$\ell\to\infty$
to~$Uu$ $\widehat\eta_\alpha$-weakly and consequently $\tilde\eta
_\alpha$-weakly.
Since $(\tilde\eta_\alpha, \dom^0)$ is a nonnegative definite
symmetric bilinear
form, it follows that
%
%
\begin{eqnarray}\label{U3}
\eta_\alpha(Uu,Uu)&=&\tilde\eta_\alpha(Uu,Uu)\le
\liminf_{\ell\to\infty}\tilde\eta_\alpha(Uu_\ell,Uu_\ell)\nonumber\\[-8pt]\\[-8pt]
&=&\liminf_{\ell\to\infty}\eta_\alpha(Uu_\ell,Uu_\ell).\nonumber
\end{eqnarray}
We can then obtain (\ref{equ1.2}) for $u\in\dom^0$ from (\ref{U1}),
(\ref
{U2}) and
(\ref{U3}) as
\begin{eqnarray*}
\eta(Uu, u-Uu)&\ge& \lim_{\ell\to\infty}\eta(Uu_\ell, u_\ell)-
\liminf_{\ell\to\infty}\eta(Uu_\ell,Uu_\ell) \\
&=& \limsup_{\ell\to\infty} \eta(Uu_\ell, u_\ell-Uu_\ell)\ge0.
\end{eqnarray*}
\upqed
\end{pf*}


For the kernel $k^*$ defined by (\ref{reverse}), we have obviously
%
%
\begin{equation}\label{equ2.17}
k^*_s(x,y)=k_s(x,y) \quad\mbox{and}\quad k^*_a(x,y)=-k_a(x,y),\qquad
x, y \in E, x\not=y.\hspace*{-35pt}
\end{equation}
Hence, if the kernel $k(x,y)$ satisfies (\ref{cond-1})--(\ref{cond-4}),
so does the kernel $k^*(x,y)$. Define $\eta^*$ as in Proposition \ref{prop2.1}
with $k^*(x,y)$
in place of $k(x,y)$. The same calculations made above for $k(x,y)$
remain valid
for $k^*(x,y)$. Note also that the domain $\dom^{0*}$ is the same as
$\dom^0$ since
the symmetric form~$\form^*$ defined by~$k^*$ is also the same as~$\form$.
Thus, we can have the following corollary.
\begin{cor}\label{cor2.1}
Assume   conditions (\ref{cond-1})--(\ref{cond-4}) hold. Then the
pair $(\eta^*,\dom^0)$ is also a regular lower bounded semi-Dirichlet
form on $L^2(E;m)$.
\end{cor}

\section{Markov property of dual semigroups}\label{sec3}

First, we fix a general lower bounded closed form $(\eta, \dom)$ on $L^2(E;m)$
satisfying (B.1)--(B.3) and make several remarks on it.
The last condition (B.3) is equivalent to

\begin{longlist}[(B.3)$'$]
\item[(B.3)$'$]
$(\tilde\eta_{\beta_0}, \dom)$ is a closed symmetric
form on $L^2(E;m)$,
\end{longlist}
where $\tilde{\eta}$ denotes the symmetrization of the form $\eta\dvtx
\tilde{\eta}(u,v)= \frac12(\eta(u,v)+\eta(v,u))$. $\eta_{\beta
_0}$ is therefore
a coercive closed form in the sense of \cite{MR92}, Definition~2.4, so
that, by~\cite{MR92}, Theorem 2.8, there exist uniquely two families of linear bounded
operators $\{G_\alpha\}_{\alpha>\beta_0},
\{\widehat G_\alpha\}_{\alpha>\beta_0}$ on $L^2(E;m)$ such that, for
$\alpha>\beta_0$, $G_\alpha(L^2(E;m))$, $\widehat G_\alpha(L^2(E;m))
\subset\dom$ and
%
%
\begin{equation}\label{equ3.1}\quad
\eta_\alpha(G_\alpha f,u)=(f,u)=\eta_\alpha(u,\widehat G_\alpha
f),\qquad
f\in L^2(E;m), u\in\dom.
\end{equation}
In particlular, $G_\alpha$ and $\widehat G_\alpha$ are mutually adjoint:
%
%
\begin{equation}\label{equ3.2}
(G_\alpha g,f)=(g, \widehat G_\alpha f),\qquad f,g\in L^2(E;m),
\alpha>\beta_0.
\end{equation}
We call $\{G_\alpha;\alpha>\beta_0\}$ (resp., $\{\widehat G_\alpha
;\alpha>\beta_0\}$)
the \textit{resolvent} (resp., \textit{dual resolvent}) associated with
$(\eta
,\dom)$.

Accordingly we see in exactly the same way as the proof of Theorem 2.8 of
\cite{MR92} that there exist strongly continuous contraction semigroups
$\{S_t;t>0\}, \{\widehat S_t; t>0\}$ of linear operators on $L^2(E;m)$
such that, for $\alpha>0, f\in L^2(E;m)$,
\[
G_{\beta_0+\alpha}f=\int_0^\infty e^{-\alpha t}S_tf\,dt,\qquad
{\widehat G}_{\beta_0+\alpha}f=\int_0^\infty e^{-\alpha t}\widehat S_tf\,dt.
\]
We then set $T_t=e^{\beta_0 t}S_t, \widehat T_t=e^{\beta_0
t}\widehat S_t$ to
get strongly continuous semigroups $\{T_t;t>0\}, \{\widehat T_t;t>0\}$
satisfying
%
%
\begin{equation}\label{equ3.3}\quad
G_\alpha f=\int_0^\infty e^{-\alpha t}T_tf\,dt,\qquad
{\widehat G}_\alpha f=\int_0^\infty e^{-\alpha t}\widehat
T_tf\,dt,\qquad
\alpha>\beta_0,
\end{equation}
as well as (\ref{equ1.1}).

We call $\{T_t;t>0\}$ (resp., $\{\widehat T_t; t>0\}$) the \textit
{semigroup} (resp.,
{dual semigroup}) on $L^2(E;m)$ associated with the lower bounded
closed form
$(\eta,\dom)$. We introduce the \textit{dual form} $\widehat\eta$ of
$\eta$ by
\[
\widehat\eta(u,v)=\eta(v,u),\qquad u,v \in\dom.
\]
Then $(\widehat\eta,\dom)$ is a lower bounded closed form on
$L^2(E;m)$ with
which $\{\widehat T_t;\break t>0\}$ and $\{\widehat G_\alpha;\alpha>\beta
_0\}$ are the
associated semigroup and resolvent, respectively.

Suppose $(\eta, \dom)$ is a lower bounded semi-Dirichlet form,
namely, it satisfies
the contraction property (\ref{equ1.2}) additionally. As in the proof of
the corollary to
Theorem 4.1 of \cite{K69} or the proof of Theorem 4.4 of \cite{MR92},
we can then
readily verify that the family $\{\alpha G_\alpha; \alpha>\beta_0\}$
is Markovian,
which is in turn equivalent to the Markovian property of $\{T_t;t>0\}$. Together
with $\{T_t;\break t>0\}$, its Laplace transform then determines a bounded
linear operator~$G_\alpha$ on $L^\infty(E;m)$ for every $\alpha>0$ and $\{\alpha
G_\alpha;
\alpha>0\}$ becomes Markovian. Further, $\{\widehat T_t;t>0\}$ is positivity
preserving in view of (\ref{equ1.1}).

Suppose additionally that $(\eta,\dom)$ is regular.
Then the associated Markovian semigroup and resolvent
can be represented by the transition function $\{P_t;t>0\}$ and the resolvent
$\{R_\alpha;\alpha>0\}$ of the associated Hunt process $X$ specified in
Theorem 2 of the next section: $P_tf=T_tf, t>0$, and $R_\alpha
f=G_\alpha f, \alpha>0$,
for any $f\in\mathcal{B}_b(E)\cap L^2(E;m)$. We call a $\sigma
$-finite measure~$\mu$ on $E$
\textit{excessive} relative to $X$ if $\mu P_t\le\mu$ for any $t>0$.
The next lemma was already observed in Silverstein~\cite{S77}.\vspace*{-1pt}
\begin{lemma}\label{lemma3.1} Let $\eta$ be a regular lower bounded
semi-Dirichlet
form on
$L^2(E;m)$.

\begin{enumerate}[(ii)]
\item[(i)]
The following three conditions are mutually equivalent:
\begin{enumerate}[1.]
\item[1.] $m$ is excessive relative to $X$.
\item[2.]  The dual semigroup $\{\widehat T_t\dvtx t>0\}$ is Markovian.
\item[3.] $\eta(u-Uu,Uu)\ge0 \mbox{ for any} u\in\dom$.
\end{enumerate}

\item[(ii)] If one of the three conditions in \textup{(i)} is satisfied,
then $\eta$
is nonnegative definite and the constant $\beta_0$ in   conditions
\textup{(B.1)}, \textup{(B.3)}
[resp., \textup{(B.2)}] can be retaken to be $0$
(resp., $1$).\vspace*{-1pt}
\end{enumerate}
\end{lemma}
\begin{pf}
(i) 3 is the Markovian criterion (\ref{equ1.2}) for the dual semigroup.
If~2 is satisfied, then for any $f\in L^2(E;m)$ with $0\le f\le1$,
$0\le\widehat
T_t f\le1$ so that $(f, P_th)=(\widehat T_tf, h)\le(1,h)$ for any
$h\in\mathcal{B}_+
\cap L^2(E;m)$, from which 1 follows. The converse can be shown
similarly.\vadjust{\goodbreak}

(ii) By the Schwarz inequality,
\[
(R_\alpha f(x))^2\le R_\alpha1(x) R_\alpha f^2(x)\le\frac1\alpha
R_\alpha f^2(x),\qquad
x\in E, f\in\mathcal{B}_b(E)\cap L^2(E;m).
\]
Assuming 1 of (i), an integration with respect to $m$ yields $\alpha
^2 \Vert G_\alpha f\Vert^2\le\Vert f\Vert^2$,
the $L^2$-contraction property of $\alpha G_\alpha$.
In view of \cite{MR92}, Theorem 2.13,
$\eta(u,u)=\lim_{\alpha\to\infty}\alpha(u-\alpha G_\alpha u, u)
u\in\dom$, which particularly implies that $\eta(u,u)\ge0,\allowbreak u\in
\dom$,
and $\{\eta_\alpha; \alpha>0\}$ become equivalent on $\dom$.
\end{pf}

We now return to the setting of the preceding section that $(\eta,
\dom^0)$ is
defined in terms of the kernel $k$ satisfying conditions
(\ref{cond-1})--(\ref{cond-4}). By Proposition \ref{prop2.1}, $\widehat\eta(u,v)=\frac12
\form(v,u)+ B(v,u)$
where $B$ is defined by (\ref{B}) on $\dom^0\times\dom^0$.
On the other hand, we have from (\ref{equ2.17}) that $\eta^*(u,v)=\frac12
\form(u,v)-B(u,v)$
and consequently
%
%
\begin{equation}\label{equ3.4}
\widehat\eta(u,v)=\eta^*(u,v)+ \bigl(B(u,v)+B(v,u)\bigr),\qquad
u,v\in\dom^0.
\end{equation}

We know from Theorem \ref{theo2.1} and Corollary \ref{cor2.1} that both
$(\eta,\dom^0)$ and
$(\eta^*,\dom^0)$ are regular lower bounded semi-Dirichlet forms.
In order to get a similar property for the dual form $\widehat\eta$,
we need to
impose on the kernel~$k$ stronger conditions than (\ref
{cond-1})--(\ref{cond-4})
making the additional term on the right-hand side of (\ref{equ3.4})
controllable.

In the rest of this section, we assume that the kernel $k$ satisfies
the condition
%
%
\begin{eqnarray}\label{cond-1-1}\quad
M_s\in L^2_{\mathrm{loc}}(E;m) \qquad\mbox{for } M_s(x)= \int_{y\not=x}
\bigl(1 \wedge d(x,y)\bigr) k_s(x,y)m(dy),\nonumber\\[-8pt]\\[-8pt]
&&\eqntext{x\in E,}
\end{eqnarray}
in place of (\ref{cond-1}), and further satisfies  condition (\ref{cond-2})
as well as (\ref{cond-3}) for $\gamma=1$ so that
%
%
\begin{eqnarray}\label{cond-2-1}
\frac{\beta_1}{2} :\!&=&\sup_{x\in E} \int_{x\not=y} | k_a(x,y)
| m(dy)\nonumber\\[-8pt]\\[-8pt]
&=&\sup_{x\in E} \frac12 \int_{x\not=y} | k(x,y)-k(y,x) |
m(dy) <\infty.\nonumber
\end{eqnarray}
Notice that  condition (\ref{cond-4}) for $\gamma=1$ is always
satisfied with
$C_3=1$.

Then the integrals
%
%
\begin{eqnarray}\label{equ3.7}
\mathcal{L}u(x) &=& \int_{y\not=x} \bigl(u(y)-u(x)\bigr) k(x,y)m(dy)
\quad\mbox{and} \nonumber\\[-8pt]\\[-8pt]
\mathcal{L}^*u(x) &=& \int_{y\not=x} \bigl(u(y)-u(x)\bigr)
k^*(x,y)m(dy),\nonumber
\end{eqnarray}
converge for $u\in C_0^{\mathrm{lip}}(E), x\in E$, and we get from
Proposition \ref{prop2.1} the
identities
%
%
\begin{equation}\label{equ3.8}\quad
\eta(u,v)=-(\mathcal{L}u, v),\qquad \eta^*(u,v)=-(\mathcal{L}^*u,
v),\qquad
u,v\in C_0^{\mathrm{lip}}(E).
\end{equation}
Furthermore,
%
%
\begin{eqnarray}\label{equ3.9}
K(x):\!&=& 2 \int_{y\not=y} k_a(x,y)m(dy)\nonumber\\[-8pt]\\[-8pt]
&=&\int_{y\not=x} \bigl(k(x,y)-k(y,x)\bigr)m(dy),\qquad x\in E,\nonumber
\end{eqnarray}
defines a bounded function on $E$ and (\ref{equ3.4}) readily leads us to
\[
\widehat\eta(u,v)=\eta^*(u,v)+ (u,Kv),\qquad
u,v\in\dom^0,
\]
which combined with (\ref{equ3.7}) means that $\widehat{\mathcal
{L}}=\mathcal
{L}^*-K$ is the formal
adjoint of~$\mathcal{L}$. $\widehat\eta$~does not necessarily
satisfy the
contraction property (\ref{equ1.2}), but the form
\[
\widehat\eta_{\beta}(u,v)=\eta^*(u,v)+\bigl(u,(K+\beta)v\bigr),\qquad \beta
\ge\beta_1,
\]
does because so does the form $\eta^*$ by Corollary \ref{cor2.1} and
$K+\beta
\ge0$ if
$\beta\ge\beta_1$. So we have the following proposition.
\begin{prop}\label{prop3.1} Assume that (\ref{cond-1-1}) and
(\ref{cond-2-1}) hold. Then $(\widehat\eta_{\beta}, \dom^0)$, which
is the dual of $(\eta_{\beta}, \dom^0)$, is a regular lower bounded
semi-Dirichlet form on $L^2(E;m)$ provided that $\beta\ge\beta_1$.
\end{prop}

This proposition means that, under conditions (\ref{cond-1-1}) and
(\ref{cond-2-1}),
$\{e^{-\beta t}\widehat T_t;\allowbreak t>0\}$ is Markovian for the dual semigroup
$\{\widehat T_t;t>0\}$ associated with $\eta$ when $\beta\ge\beta
_1$. If
(\ref{cond-2-1}) fails, the dual semigroup of $\{e^{-\beta t} T_t;
t>0\}$ may not
be Markovian no matter how large $\beta$ is.

A nonnegative Borel function $k$ on $E\times E$ is said to be a \textit
{probability
kernel} if $\int_Ek(x,y)m(dy)=1, x\in E$. A probability kernel
$k$ with the additional property
%
%
\begin{equation}\label{equ3.10}
\sup_{x\in E}\int_D k(y,x)m(dy)<\infty
\end{equation}
satisfies conditions (\ref{cond-1-1}) and (\ref{cond-2-1}) and $\eta$
defined by (\ref{equ3.8})
yields a regular lower bounded
semi-Dirichlet form on $L^2(E;m)$. We now give an example of a such a
kernel on $\mathbb R^1$
with $m$ being the Lebesgue measure for which the associated
semi-Dirichlet form $\eta$
is \textit{not} nonnegative definite so that, according to Lemma \ref
{lemma3.1}, the
associated dual semigroup
$\{\widehat T_t, t>0\}$ is \textit{not} Markovian although $\{e^{-\beta
t}\widehat T_t; t>0\}$
is Markovian for a large $\beta>0$ in view of Proposition \ref
{prop3.1}. A
transition probability density
function with respect to the Lebesgue measure of the one-dimensional Brownian
motion with a~mildly localized drift serves to be an example of such a
kernel $k$.

Consider a diffusion $Y$ on $\mathbb R^1$ with generator $\mathcal
{G}u=\frac12 u''
+\lambda b(x)u'$ where $\lambda$ is a positive constant and $b$ is a
function in
$C_0^1(\mathbb R^1)$ not identically $0$. Then
$\mathcal{G}=\frac{d}{dm}\cdot\frac{d}{ds}$ for
\[
dm(x)=m(x)\,dx,\qquad ds(x)=2 m(x)^{-1}\,dx,
\]
where
\[
m(x)=2\exp\biggl\{2\lambda\int_0^xb(y)\,dy\biggr\},
\]
namely, $Y$ is a diffusion with canonical scale $s$ and canonical
(speed) measure~$dm$.

The following facts about $Y$ are taken from \cite{F10}. Since $m(x)$
is bounded
from above and from below by positive constants, both $\pm\infty$ are
nonapproachable in the sense that $s(\pm\infty)=\pm\infty$.
Therefore, $Y$ is
recurrent and consequently conservative: $q_t(x,E)=1, x\in E$, where
$\{q_t;t>0\}$
denotes the transition function of $Y$. $Y$~is $m$-symmetric and its Dirichlet
form $(\form^Y, \dom^Y)$ on $L^2(\mathbb R^1,m)$ is given by
\[
\cases{
\displaystyle \form^Y(u,v)=\frac12\int_{\mathbb R^1}u'(x)v'(x)m(x)\,dx,\vspace*{2pt}\cr
\dom^Y=\{u\in L^2(\mathbb R^1;m)\dvtx u \mbox{ is absolutely}\vspace*{2pt}\cr
\hspace*{35.4pt}
\mbox{continuous and } \form^Y(u,u)<\infty\}\mbox{ } (\mbox{$=$}H^1(\mathbb R^1)).}
\]

For $u\in C_0^1(\mathbb R^1)$, $\form^Y(u, \frac{u}{m})$ is seen
to be equal to
$\frac12\int_{\mathbb R^1}((u')^2-2\lambda bu'u)\,dx$ and so
\[
\form^Y\biggl(u, \frac{u}{m}\biggr)= \frac12\biggl(\int_{\mathbb R^1}(u')^2\,dx+
\lambda\int_{\mathbb R^1}b'u^2\,dx\biggr).
\]
There is a finite interval $I\subset\mathbb R^1$ where $b'$ is
strictly negative.
Choose $u_0\in C_0^1(\mathbb R^1)$ not identically zero and with
support being
contained in $I$. We can then make a choice of $\lambda>0$ such that
the right-hand side of the above equation is negative for $u=u_0$.

Since $q_t$ maps $L^2(\mathbb R^1;m)$ into $\dom^Y\subset C(\mathbb R^1)$,
$q_t(x, \cdot)$ is absolutely continuous with respect to $m$ and hence with
respect to the Lebesgue measure for each \mbox{$x\in\mathbb R^1$}.
Denote by $q_t(x,y)$ its density with respect to the Lebesgue measure
so that
$\int_{\mathbb R^1}q_t(x,y)\,dy=1, x\in\mathbb R^1$, with
%
%
\begin{equation}\label{equ3.11}
q_t(y,x)=m(x)q_t(x,y)\frac{1}{m(y)}.
\end{equation}
We know that the left-hand
side of the above equation equals
\[
\lim_{t\downarrow0}\frac1t \int_{\mathbb R^1} \bigl(u(x)-q_tu(x)\bigr)\frac
{u(x)}{m(x)}m(x)\,dx
=\lim_{t\downarrow0}\frac1t \int_{\mathbb R^1}
\bigl(u(x)-q_tu(x)\bigr)u(x)\,dx\vadjust{\goodbreak}
\]
and so, for $k(x,y)=q_{t_0}(x,y)$ with a sufficiently small $t_0>0$,
\[
\eta(u_0,u_0)= -\int_{\mathbb R^1}\biggl[\int_{\mathbb
R^1}\bigl(u_0(y)-u_0(x)\bigr)k(x,y)\,dy
\biggr]u_0(x)\,dx<0.
\]
Equality (\ref{equ3.10}) follows from (\ref{equ3.11}).

\section{Associated Hunt process and martingale problem}\label{sec4}

Let $(\eta, \dom)$ be a~regular lower bounded semi-Dirichlet form on
$L^2(E;m)$ as
is defined in Section \ref{sec1}. For the symmetrization $\tilde\eta$,
$(\tilde\eta_{\beta_0},
\dom)$ is then a closed symmetric form on $L^2(E;m)$ but not
necessarily a symmetric
Dirichlet form. A symmetric Dirichlet form $\form$ on $L^2(E;m)$ with
domain $\dom$
will be called a \textit{reference} (\textit{symmetric Dirichlet})
\textit{form}
of $\eta$ if,
for each fixed $\alpha>\beta_0$,
%
%
\begin{equation}\label{equ4.1}
c_1 \form_1(u,u)\le\eta_\alpha(u,u)\le c_2 \form_1(u,u),\qquad
u\in\dom,
\end{equation}
for some positive $c_1, c_2$ independent of $u\in\dom$. $\form$ is then
a regular Dirichlet form. In what follows, we assume that $\eta$ admits
a reference form $\form$. This assumption is really unnecessary (cf.
\cite{MOR94,O88}) but convenient to simplify some arguments. The
regular lower bounded semi-Diriclet form $(\eta,\dom^0)$ constructed in
Section \ref{sec2} from a kernel $k$ satisfying (\ref{cond-1})--(\ref{cond-4})
has a reference form $(\form, \dom^0)$ defined right after~(\ref{symmetric-form}).

In formulating an association of a Hunt process with $\eta$, Carrillo Menendez
adopted a functional capacity theorem due to Ancona \cite{A72}.
More specifically, denote by $\mathcal{O}$ the family of all open sets
$A\subset E$
with $\mathcal{L}_A=\{u\in\dom\dvtx u\ge1$ $m\mbox{-a.e. on } A\}
\neq\varnothing$.
Fix $\alpha>\beta_0$ and, for $A\in\mathcal{O}$, let $e_A$ be the
$\eta_\alpha$-projection of~$0$ on $\mathcal{L}_A$ in Stampacchia's
sense \cite{S64}
(cf. \cite{MR92}, Theorem 2.6):
%
%
\begin{equation}\label{equ4.2}
e_A\in\mathcal{L}_A,\qquad \eta_\alpha(e_A,w)\ge\eta_\alpha(e_A,
e_A)\qquad
\mbox{for any } w\in\mathcal{L}_A.
\end{equation}
A set $N\subset E$ is called $\eta$-\textit{polar} if there exist decreasing
$A_n\in\mathcal{O}$ containing $N$ such that $e_{A_n}$ is $\eta
_\alpha$-convergent to
$0$ as $n\to\infty$. A numerical function $u$ on $E$ is called
$\eta$-\textit{quasi-continuous} if there exist decreasing $A_n\in
\mathcal{O}$ such that
$e_{A_n}$ is $\eta_\alpha$-convergent to $0$ as $n\to\infty$ and
$u|_{E\setminus A_n}$ is
continuous for each $n$.

The capacity Cap for the reference form $\form$ is defined by
\[
\operatorname{Cap}(A)=\inf\{\form_1(u,u)\dvtx u\in\mathcal{L}_A\},\qquad A\in
\mathcal{O}.
\]
It then follows from (\ref{equ4.1}) that
%
%
\begin{eqnarray}\label{equ4.3}
c_1 \operatorname{Cap}(A)&\le&\eta_\alpha(e_A,e_A)\le c_2 K_\alpha^2
\operatorname{Cap}(A),\qquad
A\in\mathcal{O}, \nonumber\\[-8pt]\\[-8pt]
K_\alpha&=&K+\frac{\alpha}{\alpha-\beta_0},\nonumber
\end{eqnarray}
because (\ref{equ4.2}) and (B.2) imply
$\eta_\alpha(e_A,e_A)\le K_\alpha^2\eta_\alpha(w,w), w \in
\mathcal{L}_A$.
Equa-\break tion~(\ref{equ4.3}) means that a set $N$ is $\eta$-polar iff it is
$\form$-polar in
the sense that
$\operatorname{Cap}(N)=0$, and a function $u$ is $\eta
$-quasi-continuous iff it is
$\form$-quasi-continuous in the sense that there exist decreasing
$A_n\in\mathcal{O}$
with $\operatorname{Cap}(A_n)\downarrow0$ as \mbox{$n\to\infty$} and $u
|_{E\setminus A_n}$ is continuous
for each $n$. Every element of $\dom$ admits its $\eta$-quasi-continuous
$m$-version. If $\{u_n\}\subset\dom$ is $\eta_\alpha$-convergent to
$u \in\dom$
and if each $u_n$ is $\eta$-quasi-continuous, then (\ref{equ4.1})
implies that
a subsequence
of $\{u_n\}$ converges $\eta$-q.e., namely, outside some $\eta$-polar
set, to an
$\eta$-quasi-continuous version of $u$. We shall occasionally drop
$\eta$ from the
terms $\eta$-polar, $\eta$-q.e. and $\eta$-quasi-continuity for simplicity.

Recall that the $L^2$-resolvent $\{G_\alpha;\alpha>\beta_0\}$
associated with
$\eta$ determines the resolvent $\{G_\alpha;\alpha>0\}$ on $L^\infty
(E;m)$ with
$\Vert G_\alpha f\Vert_\infty\le\frac1\alpha\Vert f\Vert_\infty$,
$\alpha>0$, $f\in L^\infty(E;m)$.
\begin{lemma}\label{lemma4.1} Suppose $G_\beta f$ admits a
quasi-continuous $m$-version
$R_\beta f$
for a fixed $\beta>\beta_0$ and for every bounded Borel $f\in
L^2(E;m)$. Then, for any~$\alpha$ with $0<\alpha\le\beta_0$ and for any bounded Borel
$f\in L^2(E;m)$,
\[
R_\alpha f(x)=\sum_{k=1}^\infty(\beta-\alpha)^{k-1}R_\beta^kf(x)
\]
converges q.e. and defines a quasi-continuous $m$-version of $G_\alpha f$.
Further the resolvent equation
\[
R_\alpha f-R_\beta f+(\alpha-\beta) R_\alpha R_\beta f=0
\]
holds q.e. for any bounded Borel $f\in L^2(E;m)$.
\end{lemma}
\begin{pf}
$\!\!$Choose a regular nest $\{F_\ell\}$ so that $R_\beta^k
f\,{\in}\,C(\{F_\ell\})$ for $k\,{\ge}\,1$. Defi\-ne $v_n(x)\,{=}\,\sum_{k=1}^n (\beta
\,{-}\,\alpha)^{k-1}
R_\beta^kf(x)$. By the resolvent equation for \mbox{$\{G_\alpha;\alpha\,{>}\,0\}$},
we have
\[
G_\alpha f= v_n+ (\beta-\alpha)^n G_\beta^nG_\alpha f.
\]
The $L^\infty$-norm of the second term of the right-hand side is
dominated by
${\frac1\alpha(\frac{\beta-\alpha}{\beta})^n \Vert f\Vert_\infty
}$, which tends to
$0$ as $n\to\infty$. Therefore, $\{v_n\}$ is convergent uniformly on
each set
$F_\ell$ to a quasi-continuous version of $G_\alpha f$. The resolvent
equation is
clear.
\end{pf}
\begin{theorem}\label{theo4.1} There exist a Borel $\eta$-polar set
$N_0\subset E$
and a Hunt process
$X=(X_t,P_x)$ on $E\setminus N_0$ which is properly associated with
$(\eta,\dom)$
in the sense that $R_\alpha f$ is a quasi continuous version of
$G_\alpha f$
for any $\alpha>0$ and any bounded Borel $f\in L^2(E;m)$. Here
$R_\alpha$ is the
resolvent of $X$ and $G_\alpha$ is the resolvent associated with $\eta$.
\end{theorem}

This theorem was proved in \cite{C75} first by assuming that $\beta
_0=0$ and then
reducing the situation to this case. Actually the proof can be carried
out without
such a reduction. Indeed, after constructing the kernel $\widetilde
V_\lambda$ of
\cite{C75}, Proposition II.2.1, for every rational $\lambda>\beta_0$
(\cite{C75}, Proposition II.2.2)
can be shown first for every rational $\lambda>\beta_0$, and then for every
$0<\lambda\le\beta_0$ by using Lemma \ref{lemma4.1}. The rest of the arguments
in \cite{C75}
then works in getting to Theorem \ref{theo4.1}.

Our next concern will be exceptional sets and fine continuity for the
Hunt process
$X=(X_t,P_x)$ appearing in Theorem \ref{theo4.1}. Denote by $\mathcal
{B}(E)$ the
family of all
Borel sets of $E$. For $B\in\mathcal{B}(E)$, we let
\[
\sigma_B=\inf\{t>0\dvtx X_t\in B\},\qquad \widehat\sigma_B=\inf\{
t>0\dvtx X_{t-}\in B\},\qquad
\inf\varnothing=\infty.
\]
$A\in\mathcal{B}(E)$ is called $X$-\textit{invariant} if
\[
P_x(\sigma_{E\setminus A}\wedge\widehat\sigma_{E\setminus A}<\infty
)=0\qquad
\forall x\in A.
\]
$N\in\mathcal{B}(E)$ is called \textit{properly exceptional} (with
respect to $X$) if
$m(N)=0$ and $E\setminus N$ is $X$-invariant.

A set $N\subset E$ is called $m$-\textit{polar} if there exists
$N_1\supset N,
N_1\in\mathcal{B}(E)$ such that $P_m(\sigma_{N_1}<\infty)=0$.
Any properly exceptional set is $m$-polar.
\begin{theorem}\label{theo4.2}
\begin{longlist}[(iii)]
\item[(i)]
For $A\in\mathcal{O}$, the function
$p_A^\alpha$ defined by
$ p_A^\alpha(x)=E_x[e^{-\alpha\sigma_A}],
x\in
E\setminus N_0$, is a quasi-continuous version of $e_A, \alpha>\beta_0$.

\item[(ii)] For any $\eta$-polar set $B$, there exists a Borel
properly exceptional
set~$N$ containing $N_0\cup B$.

\item[(iii)] If $u$ is $\eta$-quasi-continuous, then there exists a
Borel properly
exceptional set $N\supset N_0$ such that, for any $x\in E\setminus N$,
%
%
\begin{equation}\label{equ4.4}
P_x\Bigl( \lim_{t'\downarrow t}u(X_{t'})=u(X_t)\ \forall t\ge0
\mbox{ and }
\lim_{t'\uparrow t}u(X_{t'})=u(X_{t-})\ \forall
t\in(0,\zeta)\Bigr)=1,\hspace*{-28pt}
\end{equation}
where $\zeta$ is the lifetime of $X$. In particular, $u$ is finely
continuous with
respect to the restricted Hunt process $X|_{E\setminus N}$.

\item[(iv)] Any $X$-semi-polar set is $\eta$-polar.

\item[(v)] A set $N\subset E$ is $\eta$-polar if and only if $N$ is $m$-polar.
\end{longlist}
\end{theorem}
\begin{pf}
(i) A function $u\in L^2(E;m)$ is said to be $\alpha
$-excessive if
\mbox{$u\ge0$}, $\beta G_{\alpha+\beta} u\le u, \beta>0$. A function $u\in
\dom$ is
$\alpha$-excessive iff $\eta_\alpha(u, v)\ge0$ for all nonnegative
$v\in\dom$
(cf. \cite{MOR94}, Theorem 2.4). In particular, $e_A$ is $\alpha
$-excessive and
further $v=e_A\wedge p_A^\alpha$ is an $\alpha$-excessive function in
$\dom$
(cf. \cite{MOR94}, Theorem 2.6). Hence, $\eta_\alpha(v, e_A-v)\ge0$.
Since $v\in\mathcal{L}_A$, $\eta_\alpha(e_A, e_A-v)\le0$ so that
$v=e_A$ and
$e_A\le p_A^\alpha$. The converse inequality can be obtained as in the
proof of
Theorem \ref{theo6.1} below by using the optional sampling theorem for
a~supermartingale but with time parameter set being a finite set.

Since the quasi-continuous function $\beta R_{\alpha+\beta}p_A^\alpha
$ converges
to $p_A^\alpha$ as $\beta\to\infty$ pointwise and in $\eta_\alpha
$, we get the
quasi-continuity of $p_A^\alpha$.

\mbox{}\phantom{i}(ii) Choose a decreasing sets $A_n\in\mathcal{O}$ with $A_n\supset
B,
\operatorname{Cap}(A_n)\to0, n\to\infty$ and put $B_1=\bigcap_n A_n$.
By (\ref{equ4.1})
and (i),
$\lim_{n\to\infty}p_{A_n}^\alpha=0$ q.e. so that
\[
P_x(\sigma_{B_1} \wedge\widehat\sigma_{B_1}<\infty
)=0,\qquad
x\in E\setminus N_1,
\]
for some polar set $N_1$. Choose next a decreasing sets $A_n'\in
\mathcal{O}$
containing $B_1\cup N_1\cup N_0$ with $\operatorname{Cap}(A_n')\to0,
n\to\infty$
and put
$B_2=\bigcap_n A_n'$. Then the above identity holds for $x\in
E\setminus B_2$.
Moreover, the above identity holds true for $B_2$ in place of $B_1$ and
for some
polar set $N_2$ in place of $N_1$. Repeating this procedure, we get an
increasing
sequence $\{B_k\}$ of $G_\delta$-sets which are polar sets such that
\[
P_x(\sigma_{B_k}\wedge\widehat\sigma_{B_k}<\infty)=0,\qquad
x\in E\setminus B_{k+1}.
\]
It then suffices to put $N=\bigcup_k B_k$.

(iii) Choose decreasing $A_n\in\mathcal{O}$ such that Cap$(A_n)\to
0, n\to0$, and
$u|_{E\setminus A_n}$ is continuous for each $n$. Let $N$ be a properly
exceptional set constructed in~(ii) starting with this sequence $\{A_n\}$.
Then, for any $x\in E\setminus N$,\break $\lim_{n\to\infty} p_{A_n}^\alpha
(x)=0$ and
consequently $P_x(\lim_{n\to\infty}\sigma_{A_n}=\infty
)=1$, which readily
implies (\ref{equ4.4}).

(iv) We reproduce a proof by Silverstein \cite{S77}. For $B\in
\mathcal{B}(E)$,
consider the entry time $\dot\sigma_B=\inf\{t\ge0\dvtx X_t\in B\}$ and
the function
$\dot p_B^\alpha(x)=E_x[e^{-\alpha\dot\sigma_B}
],\allowbreak x\in E$,
\mbox{$\alpha>\beta_0$}. Let $K$ be a compact thin set: $K$ admits no regular
point relative
to $X$. It suffices to show that $K$ is $\eta$-polar.

Choose relatively compact open sets $\{G_n\}$ such that
$G_n\supset\overline G_{n+1}$ and $\bigcap_n G_n=K$. Due to the
quasi-left continuity
of $X$, $p_{G_n}^\alpha(x)=\dot p_{G_n}^\alpha(x)$ then
decreases to
$\dot p_K^\alpha(x)$ as $n\to\infty$ for each $x\in E$. By\vspace*{1pt} (i) and
(\ref{equ4.1}) and (\ref{equ4.2}),
the sequence $\{\dot p_{G_n}^\alpha\}$ is $\form_1$-bounded so that
the Ces\`{a}ro
mean sequence $f_n$ of its suitable subsequence is $\form_1$-convergent.
Since $f_n$ are quasi-continuous and converges to $\dot p_K^\alpha$ pointwise
as $n\to\infty$, we conclude that $\dot p_K^\alpha$ is a
quasi-continuous element
of $\dom$. On the other hand, the quasi-continuous function
$\beta R_{\alpha+\beta}\dot p_K^\alpha$ converges to $p_K^\alpha$ as
$\beta\to\infty$ pointwise and in $\eta_\alpha$ so that
$p_K^\alpha$ is also a~quasi-continuous version of $\dot p_K^\alpha$. Therefore, $p_K^\alpha
=\dot p_K^\alpha$
q.e. and in particular~$K$ is $\eta$-polar.

(v) ``only if'' part follows from (ii). To show ``if''
part, assume that
$K$ is a~compact $m$-polar set. Then $p_K^\alpha=0$ $m$-a.e. Choose for $K$
relatively compact open sets $\{G_n\}$ as in the proof of (iv) so that
the Ces\`{a}ro
mean $f_\ell$ of a certain subsequence $\{p_{G_{n_\ell}}^\alpha\}$ is
$\form_1$-convergent to $p_K^\alpha$ as $\ell\to\infty$ which is
now a zero element
of $\dom^0$. Since $f_\ell\ge1$ $m$-a.e. on $G_{n_\ell}$, we have
Cap$(K)\le$
Cap$(G_{n_\ell})\le\form_1(f_\ell,f_\ell)$ and we get Cap$(K)=0$
by letting
$\ell\to\infty$. For any Borel $m$-polar set $N$, we have
Cap$(N)=\sup\{\operatorname{Cap}(K)\dvtx K\subset N, K \mbox{is
compact}\}=0$.
\end{pf}

Clearly, the restriction of $X$ outside its properly exceptional set is
again a Hunt
process properly associated with $\eta$.

Our final task in this section is to relate the Hunt process of Theorem
\ref{theo4.1} to a
martingale problem.\vadjust{\goodbreak}

We consider the case where $\eta$ admits the expression
%
%
\begin{equation}\label{equ4.5}
\eta(f,g)=-(\el f, g),\qquad f\in\D(\el), g\in\dom,
\end{equation}
for a operator $\el$ with domain $\D(\el)$ satisfying the
following:\vspace*{8pt}

(L.1) $\D(\el)$ is a linear subspace of $\dom\cap C_0(E)$,

(L.2) $\el$ is a linear operator sending $\D(\el)$ into $L^2(E;m)
\cap C_b(E)$,

(L.3) there exists a countable subfamily $\D_0$ of $\D(\el)$
such that
each $f\in\D(\el)$ admits $f_n\in\D_0$ such that $f_n, \el f_n$
are uniformly
bounded and converge pointwise to $f, \el f$, respectively, as $n\to
\infty$.

We also consider an additional condition that

(L.4) there exists $f_n\in\D(\el)$ such that $f_n, \el f_n$
are uniformly
bounded and converge to $1, 0$, respectively, as $n\to\infty$.
\begin{theorem}\label{theo4.3}
Assume that $\eta$ admits the expression (\ref{equ4.5})
with $\el$
satisfying conditions \textup{(L.1)}, \textup{(L.2)}, \textup{(L.3)}.

\begin{longlist}[(ii)]
\item[(i)]
There exists then a Borel properly exceptional set $N$
containing $N_0$ such that,
for every $f\in\D(\el)$,
%
%
\begin{equation}\label{equ4.6}
M_t^{[f]}=f(X_t)-f(X_0)-\int_0^t (\el f)(X_s)\,ds,\qquad t\ge0,
\end{equation}
is a $P_x$-martingale for each $x\in E\setminus N$.

\item[(ii)] If the additional condition \textup{(L.4)} is satisfied,
then the Hunt
process $X|_{E\setminus N}$ is conservative.
\end{longlist}
\end{theorem}
\begin{pf}
(i) Take $f\in\D(\el)$ and $g\in L^2(E;m)$. By (\ref{equ4.5})
and (\ref{equ3.2}),
we have, for $\alpha>\beta_0$,
\begin{eqnarray*}
(G_{\alpha}\mathcal{L} f, g) &=& (\mathcal{L}f, \widehat G_\alpha g)
=-\eta(f, \widehat G_{\alpha}g)\\
&=& -\eta_{\alpha}(f, \widehat G_{\alpha}g)+\alpha(f, \widehat
G_{\alpha}g)\\
&=&-(f,g)+\alpha(G_{\alpha}f, g).
\end{eqnarray*}
Thus,
$ {(G_{\alpha}\mathcal{L}f, g)=(\alpha G_{\alpha}f- f, g)}$
holds for any
$g \in\dom$ and
\[
\frac1{\alpha} G_{\alpha}(\mathcal{L}f)(x)=
G_{\alpha}f(x)- \frac{f(x)}{\alpha},\qquad m\mbox{-a.e.}
\]
We denote by $\{P_t;t\ge0\}$ and $\{R_\alpha;\alpha>0\}$ the
transition function
and the resolvent of $X$, respectively:
\[
P_th(x)={\mathbb E}_x[h(X_t)],\qquad R_{\alpha}h(x)=\int_0^{\infty}
e^{-\alpha t} P_th(x)\,dt.
\]
Since $X$ is properly associated with $\eta$ by Theorem \ref{theo4.1},
we get
\[
\frac1{\alpha} R_{\alpha}(\mathcal{L}f)(x)=
R_{\alpha}f(x)- \frac{f(x)}{\alpha}, \qquad\mbox{q.e.}
\]
Hence, by virtue of Theorem \ref{theo4.2}(ii), there exists a Borel properly
exceptional set
$N$ such that
\[
\int_0^{\infty}e^{-\alpha t} \biggl(\int_0^t P_s (\mathcal
{L}f)(x)\,ds\biggr)\,dt
=\int_0^{\infty} e^{-\alpha t} \bigl(P_tf(x)-f(x)\bigr)\,dt,\qquad x\in
E\setminus N,
\]
holds for any $\alpha\in{\mathbb Q}_+$ with $\alpha>\beta_0$ and
for any
$f\in\D_0$.

Since $P_t h(x)$ is a right continuous in $t\ge0$ for any $h\in
C_b(E)$, we get
%
%
\begin{equation}
\label{equ4.7}
P_tf(x)-f(x)=\int_0^t P_s(\mathcal{L}f)(x)\,ds,\qquad t\ge0, x\in
E\setminus N,
\end{equation}
holding for any $f\in\D_0$. By virtue of  condition (L.3),
we conclude
that the equation (\ref{equ4.7}) holds true for any $f\in\D(\el)$.
Equation (\ref{equ4.7})
implies that,
for any $f\in\D(\el)$, the functional $M^{[f]}_t, t\ge0$, defined
by (\ref{equ4.6})
is a mean zero, square integrable additive functional of the Hunt process
$X|_{E\setminus N}$ so that it is a $P_x$-martingale for each
$x\in E\setminus N$.

(ii) Under the additional condition (L.4), we let $n\to\infty$
in   equation
(\ref{equ4.7}) with $f_n$ in place of $f$ arriving at $P_t1=1, t\ge0$.
\end{pf}

Theorem \ref{theo4.3} will enable us in the next section to relate our Hunt
process to the
solution of a martingale problem in a specific case.

\section{Stable-like process}\label{sec5}

In this section, we consider the case that $E=\rd$ and $m(dx)=dx$ is the
Lebesgue measure on $\rd$. For a positive measurable function $\alpha(x)$
defined on $\rd$, Bass introduced the following integro-differential operator
in \cite{B88b} (see also \cite{B88a,B04}): for $u\in C_b^2(\rd)$,
\begin{eqnarray}
\mathcal{L}u(x)= w(x)\int_{h\not=0}\bigl(u(x+h)-u(x)-\nabla u(x)
\cdot h
\mathbf{1}_{B(1)}(h) \bigr) |h|^{-d-\alpha(x)} \,dh,\nonumber\\
&&\eqntext{x\in\rd,}
\end{eqnarray}
where $w(x)$ is a function chosen so that $\mathcal
{L}e^{iux}=-|u|^{\alpha(x)}e^{iux}$
and $C_b^2(\rd)$ denotes the set of twicely differentiable bounded functions.
If $\alpha$ is Lipschitz continuous, bounded below by a constant which
is greater
than $0$, and bounded above by a constant which is less than $2$, then he
constructed a~unique strong Markov process associated with $\mathcal
{L}$ by solving
the $\mathcal{L}$-martingale problem for every starting point $x\in
\rd$. Using the
theory of stochastic differential equation with jumps, Tsuchiya \cite
{T92} also
succeeded in constructing the Markov process associated with $\mathcal
{L}$ (see also~\cite{N94}). Note that the weight function $w(x)$ is given by
%
%
\begin{equation}\label{weight}\qquad
w(x)=\frac{{\Gamma(({1+\alpha(x)})/{2})
\Gamma(({\alpha(x)+d})/{2}) \sin({\pi
\alpha(x)}/{2})}}
{{2^{1-\alpha(x)}\pi^{d/2+1}}},\qquad x\in\rd
\end{equation}
(see, e.g., \cite{AS61}).

Put $k(x,y)=w(x)|x-y|^{-d-\alpha(x)}, x,y \in\rd$ with $x\not=y$.
Then this falls into our case when we consider the following conditions:
there exist positive constants
$\underline{\alpha}, \overline{\alpha}, M$ and $\delta$ so that
for $x,y \in\rd$,
%
%
\begin{eqnarray}\label{stable-like}
&\displaystyle 0<\underline{\alpha} \le\alpha(x) \le\overline{\alpha}<2,
\overline{\alpha} <1+\frac{\underline{\alpha}}{2} \quad\mbox{and}&
\nonumber\\[-8pt]\\[-8pt]
&\displaystyle |\alpha(x)-\alpha(y)| \le M |x-y|^{\delta}
\qquad\mbox{for } \delta\mbox{ with }
0<\frac12(2\overline{\alpha}-\underline{\alpha})
<\delta\le1.&\nonumber
\end{eqnarray}
\begin{prop}\label{prop5.1} Assume (\ref{stable-like}) holds. Then
conditions
(\ref{cond-1})--(\ref{cond-4}) are satisfied by the function
%
%
\begin{equation}\label{equ5.3}
k(x,y)=w(x)|x-y|^{-d-\alpha(x)},\qquad x,y \in\rd, x\not=y.
\end{equation}
\end{prop}
\begin{pf}
Note first that, from   equation (\ref{weight}) defining
the weight $w(x)$, we easily see that there exist constants $c_i$ $(i=1,2,3)$
so that for $x,y \in\rd$,
\[
c_1 \le w(x) \le c_2,\qquad
|w(x)-w(y)|\le c_3 |\alpha(x)-\alpha(y)|.
\]
Then
\begin{eqnarray*}
k_s(x,y) &=& \tfrac12 \bigl( w(x)|x-y|^{-d-\alpha(x)} +w(y)
|x-y|^{-d-\alpha(y)}\bigr) \\
&\le& \cases{
M |x-y|^{-d-\overline{\alpha}}, &\quad$|x-y|\le1$, \cr
M |x-y|^{-d-\underline{\alpha}}, &\quad$|x-y|>1$.}
\end{eqnarray*}
This and the condition $0<\underline{\alpha} \le\overline{\alpha}<2$
imply that condition (\ref{cond-1}) is fulfilled because the function
$M_s$ in it is bounded. Condition (\ref{cond-2}) is also valid as
$|k_a(x,y)|\le k_s(x,y)$.

On the other hand, since
\begin{eqnarray*}
k_a(x,y) &=& w(x) |x-y|^{-d-\alpha(x)} -w(y)|x-y|^{-d-\alpha
(y)}\\
&=& \bigl(w(x)-w(y)\bigr) |x-y|^{-d-\alpha(x)}\\
&&{}+w(y) |x-y|^{-d}\bigl(|x-y|^{-\alpha(x)} -|x-y|^{ -\alpha(y)}\bigr)
\end{eqnarray*}
and
\[
|x-y|^{-\alpha(x)}-|x-y|^{-\alpha(y)} =\int^{\alpha(x)}_{\alpha
(y)} |x-y|^{-u}
\frac{1}{{\ln}|x-y|^{-1}}\,du,
\]
we see that for $|x-y|< 1$,
\begin{eqnarray*}
|k_a(x,y)| &\le& |w(x)-w(y)|\cdot|x-y|^{-d-\alpha
(x)} \\
&&{} + w(y) |x-y|^{-d} |\alpha(x)-\alpha(y)| \cdot
|x-y|^{-(\alpha(x) \vee\alpha(y))} \frac{1}{{\ln}|x-y|^{-1}} \\
&\le& M \biggl(|x-y|^{-d-\overline{\alpha}+\delta} +
|x-y|^{-d-\overline{\alpha}+\delta} \frac{1}{{\ln}|x-y|^{-1}}
\biggr) \\
&\le& M' |x-y|^{-d-\overline{\alpha}+\delta} \frac{1}{{\ln}
|x-y|^{-1}}.
\end{eqnarray*}
So if $\gamma$ satisfies
\[
\gamma(d+\overline{\alpha}-\delta)-(d-1)<1,
\]
then   condition (\ref{cond-3}) holds. As for   condition (\ref
{cond-4}), note that
\[
k_s(x,y) \ge M' |x-y|^{-d-\underline{\alpha}},\qquad |x-y|<1.
\]
So, (\ref{cond-4}) is valid when
\[
(d+\overline{\alpha}-\delta)(2-\gamma) <d+\underline{\alpha}.
\]
Therefore,   conditions (\ref{cond-3}) and (\ref{cond-4}) hold
provided that
$\gamma$ satisfies
\[
\frac{d+2\overline{\alpha}-2\delta-\underline{\alpha
}}{d+\overline{\alpha}-\delta}
< \gamma<\frac{d}{d+\overline{\alpha}-\delta}.
\]
\upqed
\end{pf}

Let $(\eta,\dom^0)$ be the regular lower bounded semi-Dirichlet form
on $L^2(\rd)$
associated with the kernel (\ref{equ5.3}) satisfying (\ref
{stable-like}) according to Theorem \ref{theo2.1}.
Let $X=(X_t,P_x)$ be the Hunt process on $\rd$ properly associated
with $(\eta,
\dom)$ by Theorem \ref{theo4.1}.

Define a linear operator $\el$ by
%
%
\begin{equation}\label{equ5.4}
\cases{
\D(\el)=C_0^2(\rd),\vspace*{2pt}\cr
\displaystyle {\mathcal L}u(x)=\int_{h\not=0} \bigl(u(x+h)-u(x)-\nabla
u(x)\cdot h
\mathbf{1}_{B_1(0)}(h)\bigr)\frac{w(x) \,dh}{|h|^{d+\alpha(x)}},\vspace*{2pt}\cr
\qquad x\in
\rd.}\hspace*{-35pt}
\end{equation}
$C_0^2(\rd)$ is a linear subspace of $\dom^0\cap C_0(\rd)$ and, by
condition (\ref{stable-like}), we can see that $\el$ maps $C_0^2(\rd)$
into $L^2(\rd)\cap C_b(\rd )$. As any continuously differentiable
function and its derivatives can be simultaneously 
approximated by
polynomials and their derivatives uniformly on each
rectangles (cf. \cite{CH53}, Chapter II), conditions (L.1), (L.2),
(L.3) in the preceding section on $\el$ are fulfilled. We can easily
verify that the present $\el$ satisfies   condition (L.4) as well.

Since the vector valued function $h w(x)\mathbf{1}_{B_1(0)}(h)
|h|^{-d-\alpha(x)}$ is
odd with respect to the variable $h$ for each $x \in\rd$, we get
for $u \in C_0^2(\rd)$,
\begin{eqnarray*}
\eta^n(u,v)&=& -\iint_{|x-y|>1/n} \bigl(u(y)-u(x)\bigr)v(x)
\frac{w(x)}{|x-y|^{d+\alpha(x)}}\,dx\,dy \\
&=& -\iint_{|h|>1/n} \bigl(u(x+h)-u(x)\bigr)v(x)
\frac{w(x)}{|h|^{d+\alpha(x)}}\,dx\,dh \\
&=& -\iint_{|h|>1/n} \bigl(u(x+h)-u(x) -\nabla u(x) \cdot h
\mathbf{1}_{B_1(0)}(h)\bigr) v(x) \\
&&\hspace*{51.4pt}{}\times\frac{w(x)}{|h|^{d+\alpha(x)}}\,dx\,dh.
\end{eqnarray*}
By letting $n\to\infty$, we have
\[
\eta(u,v)=-(\mathcal{L}u,v),
\]
that is, $\eta$ is related to $\el$ by (\ref{equ4.5}).

By virtue of Theorem \ref{theo4.3}, there exists a Borel properly
exceptional set
$N \subset\rd$
so that $X|_{{\mathbb R}^d\setminus N}$ is conservative and,
for each $x\in\rd\setminus N$,
\[
M^{[f]}_t=f(X_t)-f(X_0)-\int_0^t(\el f)(X_s)\,ds,\qquad t\ge0,
\]
is a martingale under $P_x$ for every $f\in C_0^2(\rd)$. Approximating
$f\in
C_b^2(\rd)$ by a uniformly bounded sequence $\{f_n\}\subset C_0^2(\rd
)$ such that
$\{\el f_n\}$ is uniformly bounded and convergent to $\el f$, we see
that (\ref{equ4.6})
remains valid for $f\in C_b^2(\rd)$ and $M^{[f]}_t$ is still a
martingale under
${\mathbb P}_x$ for $x\in\rd\setminus N$. For each $x\in\rd
\setminus N$,
the measure ${\mathbb P}_x$ is thus a solution to the martingale
problem for the
operator $\el$ of (\ref{equ5.4}) starting at $x$ so that ${\mathbb P}_x$
coincides with the
law constructed by Bass \cite{B88b} because of the uniqueness
also due to \cite{B88b}.
\begin{remark}\label{remark5.1}
Let
%
%
\begin{equation}\label{equ5.5}
k^*(x,y)=\frac{w(y)}{|x-y|^{d+\alpha(y)}},\qquad x,\in\rd,
x\neq y.
\end{equation}
Under   condition (\ref{stable-like}), the form $\eta^*$
corresponding to the
kernel $k^*$ is a
regular lower bounded semi-Dirichlet form on $L^2(\rd)$ by virtue of
Proposition \ref{prop5.1}
and Corollary \ref{cor2.1}. By Theorem \ref{theo4.1}, $\eta^*$ admits a properly
associated Hunt process
$X^*$ on $\rd$. Furthermore, we can have 
an explicit
expression $\eta^*(u,v)=-(\el^* u,v)$ for $u\in C_0^2(\rd)$ and
$v\in\dom^0$ with
\begin{eqnarray*}
{\mathcal L}^*u(x)&=& \int_{h\not=0} \bigl(u(x+h)-u(x)-\nabla
u(x)\cdot h
\mathbf{1}_{B_1(0)}(h)\bigr)\frac{w(x+h)\,dh}{|h|^{d+\alpha(x+h)}}\\
&&{} +\frac12\int_{0<|h|<1}\nabla u(x)\cdot h
\biggl(\frac{w(x+h)}{|h|^{d+\alpha(x+h)}}-\frac{w(x-h)}{|h|^{d+\alpha
(x-h)}}\biggr)\,dh,\qquad
x\in\rd.
\end{eqnarray*}
In a lower order case as is considered in Section \ref{sec3}, both $\el
$ and
$\el^*$ admit simpler expressions (\ref{equ3.7}) and $\el^*-K$ is a formal
adjoint of
$\el$ for a function~$K$ defined by (\ref{equ3.9}).
\end{remark}

\section{Associated Hunt processes on open subsets and on their
closures}\label{sec6}
We make the same assumptions on $E, m, k$ as in Section \ref{sec2}. Let
$D$ be
an arbitrary
open subset of $E$ and $\overline D$ be the closure of $D$, $m_D$ is
defined to
be $m_D(B)=m(B\cap D), B\in\mathcal{B}(E)$ and $(u,v)_D$ denotes the
inner product\vadjust{\goodbreak}
of $L^2(D,m_D)$ $(\mbox{$=$}L^2(\overline D, m_D))$. Consider the related function spaces
$C_0^{\mathrm{lip}}(\overline D)$ and $C_0^{\mathrm{lip}}(D)$
introduced in
Section \ref{sec1}.
Define
%
%
\begin{equation}\label{equ6.1}
\cases{
\displaystyle \form_D(u,v):= \iint_{D\times D\setminus\di} \bigl(u(y)-u(x)\bigr)
\bigl(v(y)-v(x)\bigr)\vspace*{2pt}\cr
\hspace*{106pt}{}\times k_s(x,y)m_D(dx)m_D(dy),\vspace*{2pt}\cr
\displaystyle \dom^r_D=\{ u\in L^2(D;m_D)\dvtx u \mbox{ is Borel
measurable and }
\form_D(u,u)<\infty\},}\hspace*{-35pt}
\end{equation}
and let $\dom_{\bar D}$ and $\dom_D^0$ be the $\form_{D,1}$-closures of
$C_0^{\mathrm{lip}}(\overline D)$ and $C_0^{\mathrm{lip}}(D)$ in $\dom_D^r$,
respectively.
$(\form_D,\dom_{\bar D})$ [resp., $(\form_D^0, \dom_D^0)$] is a
regular symmetric
Dirichlet form on $L^2(\overline D;m_D)$ [resp., $L^2(D; m_D)$] where
$\form_D^0$
denotes the restriction of $\form_D$ to $\dom_D^0\times\dom_D^0$.
Furthermore,
in view of \cite{FOT94}, Theorem 4.4.3, we have the
identity\looseness=1
%
%
\begin{equation}\label{equ6.2}
\dom_D^0= \{u\in\dom_{\bar D}\dvtx\tilde u=0, \form
_D\mbox{-q.e. on }
\partial D \},
\end{equation}\looseness=0
where $\tilde u$ denotes an $\form_D$-quasi continuous version of
$u\in\dom_{\bar D}$. We keep in mind that a subset of $D$ is polar for
$(\form_D, \dom_D^0)$ iff so it is for $(\form_D, \dom_{\bar D})$,
and the restriction to $D$ of a quasi continuous function with respect
to the
latter is quasi-continuous with respect to the former.

Now define for $u\in C_0^{\mathrm{lip}}(\overline D)$ and $n\in{\mathbb N}$
%
%
\begin{equation}\label{equ6.3}
\mathcal{L}^n_Du(x):=\int_{\{y\in D\dvtx d(x,y)>1/n\}}
\bigl(u(y)-u(x)\bigr)k(x,y)m_D(dy),\qquad
x\in D.\hspace*{-35pt}
\end{equation}
Then, just as in Proposition \ref{prop2.1} and Theorem \ref{theo2.1} of
Section \ref{sec2}, we
conclude that the finite
limit
%
%
\begin{equation}\label{equ6.4}\qquad
\eta_D(u,v)=-\lim_{n\to\infty}\int_D \mathcal{L}_D^nu(x) v(x)
m_D(dx)\qquad\mbox{for }
u,v\in C_0^{\mathrm{lip}}(\overline D)
\end{equation}
exists, $\eta_D$ extends to $\dom_{\bar D}\times\dom_{\bar D}$ and
$(\eta_D,\dom_{\bar D})$ becomes a regular lower boun\-ded
semi-Dirichlet form on
$L^2(\overline D;m_D)$ possessing $(\form_D, \dom_{\bar D})$ as its reference
symmetric Dirichlet form. In parallel with $(\eta_D,\dom_{\bar D})$,
the space
$(\eta_D^0,\dom_D^0)$ becomes a regular lower bounded semi-Dirichlet
form on
$L^2(D;m_D)$ possessing $(\form_D^0, \dom_D^0)$ as its reference symmetric
Dirichlet form. Here $\eta_D^0$ is the restriction of $\eta_D$ to
$\dom_D^0\times\dom_D^0$.

Let $X^{\bar D}=(X_t, P_x)$ be a Hunt process on $\overline D$ properly
associated
with the form $(\eta_D,\dom_{\bar D})$ on $L^2(\overline D;m_D)$.
Denote by
$X^{D,0}=(X^{D,0}_t, P_x)$ the part process of $X^{\bar D}$ on $D$, namely,
$X^{D,0}_t$ is obtained from $X_t$ by killing upon hitting the boundary
$\partial D$:
\[
X^{D,0}_t=X_t,\qquad t<\sigma_{\partial D};\qquad
X^{D,0}_t=\Delta,\qquad
t \ge\sigma_{\partial D},
\]
$X^{D,0}$ is a Hunt process with state space $D$.
\begin{theorem}\label{theo6.1}
The part process $X^{D,0}$ of $X^{\bar D}$ on $D$ is properly
associated with the regular lower bounded semi-Dirichlet form
$(\eta^0_D, \dom_D^0)$ on $L^2(D; m_D)$.\looseness=1\vadjust{\goodbreak}
\end{theorem}

\begin{pf}
Let $\{R_\alpha; \alpha>0\}$ be the resolvent of
$X^{\bar D}$.
$\sigma$ will denote the hitting time of $\partial D$ by $X^{\bar
D}\dvtx\sigma=\sigma_{\partial D}$. Put, for $\alpha>0$ and $x\in
\overline D$,
\begin{eqnarray*}
R_\alpha^{D,0}f(x)&=&E_x\biggl[\int_0^\sigma e^{-\alpha
t}f(X_t)\,dt\biggr],\\
H^{\partial D}_\alpha u(x)&=&E_x[e^{-\alpha\sigma}u(X_\sigma
)],\qquad
x\in\overline D.
\end{eqnarray*}
$\{R^{D,0}_\alpha|_D; \alpha>0\}$ is the resolvent of the part process
$X^{D,0}$ of $X^{\bar D}$ on $D$.

We need to prove that, for any $\alpha>\beta_0$ and any $f\in
\mathcal{B}(\overline D)
\cap L^2(\overline D,m_D)$,
%
%
\begin{eqnarray}\label{equ6.5}
&&R_\alpha^{D,0}f \mbox{ is }
\eta_D^0\mbox{-quasi-continuous},\nonumber\\[-8pt]\\[-8pt]
&&R_\alpha^{D,0} f\in\dom_D^0,\qquad \eta_{D,\alpha}^0(R_\alpha^{D,0}
f,v)=(f,v)_D\qquad
\mbox{for any } v\in\dom_D^0.\nonumber\hspace*{-25pt}
\end{eqnarray}

We denote by $\mathcal G$ the space appearing in the right-hand side of
(\ref{equ6.2}).
Notice that $\form_D$-q.e. (resp., $\form_D$-quasi-continuity) is now
a synonym of
$\eta_D$-q.e. (resp., $\eta_D$-quasi-continuity). As the set of points of
$\partial D$ that are irregular for $\partial D$ is known to be semi-polar,
we have $P_x(\sigma=0)=1$ and so $R_\alpha^{D,0}f(x)=0$ for $\eta_D$-q.e.
$x\in\partial D$ owing to Theorem \ref{theo4.2}(iv). Since
\begin{eqnarray*}
&&R_\alpha f \mbox{ is } \eta_D\mbox{-quasi-continuous},\\
&&R_\alpha f\in\dom_{\bar D},\qquad \eta_{D,\alpha}(R_\alpha f,
v)=(f,v)_D\qquad
\mbox{for any } v\in\dom_{\bar D}
\end{eqnarray*}
and
%
%
\begin{equation}\label{equ6.6}
R_\alpha f(x)=R_\alpha^{D,0}f(x)+H_\alpha^{\partial D} R_\alpha
f(x),\qquad
x\in\overline D,
\end{equation}
we see that, for the proof of (\ref{equ6.5}), it is enough to show that
%
%
\begin{eqnarray}\label{equ6.7}
&&H_\alpha^{\partial D} R_\alpha f \mbox{ is }
\eta_D\mbox{-quasi-continuous},\nonumber\\[-8pt]\\[-8pt]
&&H_\alpha^{\partial D} R_\alpha f\in\dom_{\bar D},\qquad
\eta_{D,\alpha}(H^{\partial D}_\alpha R_\alpha f, v)=0 \qquad\mbox{for any }
v\in\mathcal G.\nonumber\hspace*{-25pt}
\end{eqnarray}

To this end, we fix $\alpha>\beta_0,
f\in\mathcal{B}_+(\overline D) \cap L^2(\overline D;m_D)$ and put
$u=R_\alpha f$.
Consider a closed convex subset of $\dom_{\bar D}$ defined by
\[
\mathcal{L}_{u,\partial D}=\{v\in\dom_{\bar D}, \tilde v\ge\tilde
u\mbox{ q.e. on } \partial D\}.
\]
Let $u_\alpha$ be the $\eta_{D,\alpha}$-projection of $0$ on
$\mathcal{L}_{u,\partial D}$:
\[
u_\alpha\in\mathcal{L}_{u,\partial D},\qquad \eta_{D,\alpha
}(u_\alpha, v-u_\alpha)\ge0,\qquad
\mbox{for any } v\in\mathcal{L}_{u,\partial D}.
\]
Both $u$ and $u_\alpha$ are $\alpha$-excessive elements of $\dom
_{\bar D}$.
By making use of the function $v=u_\alpha\wedge u$ as in the proof of
Proposition \ref{prop3.1}(i), we readily get
%
%
\begin{equation}\label{equ6.8}
\tilde u_\alpha= u\mbox{ q.e. on } \partial D,\qquad
\eta_{D,\alpha}(u_\alpha, v)=0\qquad \mbox{for any } v\in\mathcal{G}.
\end{equation}

Finally, we prove that
%
%
\begin{equation}\label{equ6.9}
H_\alpha^{\partial D} u \mbox{ is } \eta_D\mbox{-quasi
continuous},\qquad
H^{\partial D}_\alpha u=u_\alpha,
\end{equation}
which leads us to the desired property (\ref{equ6.7}). By (\ref
{equ6.6}), $H^{\partial
D}_\alpha u$
is an $\alpha$-excessive function dominated by $u\in\dom_{\bar D}$
so that
$H^{\partial D}_\alpha u$ is a quasi-continuous element of $\dom_{\bar D}$.
Further $H_\alpha^{\partial D} u=u$ q.e. on $\partial D$ by (\ref
{equ6.6}) and
an observation
made preceding it. Let $v=H_\alpha^{\partial D} u\wedge u_\alpha$.
Then $\tilde v=H_\alpha^{\partial D} u\wedge\tilde u_\alpha=u$ q.e.
on $\partial D$ so that
$\eta_{D,\alpha}(u_\alpha, u_\alpha-v)=0$ by (\ref{equ6.8}). On the
other hand,
$v$ is $\alpha$-excessive and so $\eta_{D,\alpha}(v, u_\alpha-v)\ge0$.
Consequently, $\eta_\alpha(u_\alpha-v, u_\alpha-v)\le0$ and we get
the inequality
$u_\alpha\le H_\alpha^{\partial D} u$.

To get the converse inequality, consider a bounded nonnegative Borel
function $h$
on $D$ with $\int_D h\,dm=1$. Denote by $\{p_t;t\ge0\}$ the
transition function
of $X^{\bar D}$. We choose a Borel measurable quasi-continuous version~%
$\tilde u_\alpha$ of $u_\alpha\in\dom_{\bar D}$. We set $\tilde
u_\alpha(\Delta)=0$
for the cemetery $\Delta$ of $X^{\bar D}$. Since $u_\alpha$ is
$\alpha$-excessive,
$e^{-\alpha t} p_t\tilde u_\alpha\le\tilde u_\alpha$ \mbox{$m$-a.e.},
and we can see
that the process $\{Y_t=e^{-\alpha t}\tilde u_\alpha(X_t); t\ge
0\}$ is a
right continuous positive supermartingale under $P_{h\cdot m}$ in view of
Theorem \ref{theo4.2}(iii). For any compact set $K\subset\partial D$,
we get
from the
optional sampling theorem and (\ref{equ6.8}),
\begin{eqnarray*}
E_{h\cdot m}[Y_{\sigma_K}]&= &
E_{h\cdot m}[e^{-\alpha\sigma_K}\tilde u_\alpha(X_{\sigma
_K})]\\
&=&E_{h\cdot m}[e^{-\alpha\sigma_K} u(X_{\sigma_K})]
\le E_{h\cdot m}[Y_0]\\
&=&(h, u_\alpha)_D.
\end{eqnarray*}
By choosing $K$ such that $\sigma_K\downarrow\sigma$ $P_{h\cdot m}$-a.e.,
we obtain $(h, H_\alpha^{\partial D} u)_D\le(h,u_\alpha)_D$ and
$H_\alpha^{\partial D} u\le u_\alpha$.
\end{pf}

As a preparation for the next lemma, we take any open set $G\subset D$
and denote by
$m_G$ the restriction of $m$ to $G$. Let $\dom_G^0$ be the $\form
_{D,1}$-closure of~$C_0^{\mathrm{lip}}(G)$ in $\dom_D^r$ and $\eta_G^0$ be the restriction of
$\eta_D$ to
$\dom_G^0\times\dom_G^0$. Then, just as above,
\[
\dom_G^0=\{u\in\dom_{\bar D}\dvtx\tilde u=0\ \form_D
\mbox{ q.e. on }
\overline D\setminus G\}
\]
and $(\eta_G^0, \dom_G^0)$ becomes a regular lower bounded
semi-Dirichlet form on
$L^2(G$; $m_G)$ with which the part process $X^{G,0}$ of $X^{\bar D}$ on
$G$ is properly
associated. The resolvent of $X^{G,0}$ will be denoted by $R_\alpha^{G,0}$.

Define
\[
H_\alpha^{\bar D\setminus G}u(x)=E_x[e^{-\alpha\sigma_{\bar
D\setminus G}} u
(X_{\sigma_{\bar D\setminus G}})],\qquad x\in
\overline D.
\]
As (\ref{equ6.7}), we have, for $u=R_\alpha f,
f\in\mathcal{B}(\overline D) \cap L^2(\overline D;m_D), \alpha>\beta_0$,
%
%
\begin{eqnarray}\label{equ6.10}
&&H_\alpha^{\bar D\setminus G} u \mbox{ is }
\eta_D\mbox{-quasi-continuous,}\nonumber\\[-8pt]\\[-8pt]
&&H_\alpha^{\bar D\setminus G} u\in\dom_{\bar D},\qquad
\eta_{D,\alpha}(H_\alpha^{\bar D\setminus G} u,v)=0 \qquad\mbox{for any }
v\in\dom_G^0,\nonumber\hspace*{-25pt}
\end{eqnarray}
and the bound $\eta_{D,\alpha}(H_\alpha^{\bar D\setminus G}u,
H_\alpha^{\bar D\setminus G}u)\le\eta_{D,\alpha}(u,u)$.
We can easily see\break that~(\ref{equ6.10}) holds true for any
$u\in\dom^{\bar D}\cap C_0(\overline D)$ where $C_0(\overline D)$
denotes the
restrictions to $\overline D$ of functions in $C_0(E)$.
In fact, by the resolvent equation,~(\ref{equ6.10}) is true for
$R_\beta u$, \mbox{$\beta>\beta_0$},
in place of $u$. Since $\{\beta_nR_{\beta_n}u\}$ converges to $u$
pointwise as well
as in $\eta_{D,\alpha}$-metric as $\beta_n\to\infty$, so does the sequence
$\{\beta_nH_\alpha^{\bar D\setminus G}R_{\beta_n}u\}$, arriving at
the validity of
(\ref{equ6.10}) for such $u$.
\begin{lemma}\label{lemma6.1}
Let $G$ be a relatively compact open set with $\overline
G\subset D$.
Then for any $v\in\dom^{\bar D}\cap C_0(\overline D)$ with
$\operatorname{supp}[v] \subset\overline D\setminus\overline G$, it
follows for
$\alpha>\beta_0$ that
%
%
\begin{equation}\label{equ6.11}
E_x[e^{-\alpha\tau_G} v(X_{\tau_G})]=
R_\alpha^{G,0}g_v(x)\qquad \mbox{for q.e. } x\in G,
\end{equation}
where $\tau_G=\sigma_{\bar D\setminus G}\wedge\zeta$ is the first
leaving time
from $G$ and $g_v$ is a function given by
%
%
\begin{equation}\label{equ6.12}
g_v(x)=1_G(x) \int_{\overline D\setminus\bar G} k(x,y)
v(y)m_D(dy),  \qquad x\in\overline D.
\end{equation}
\end{lemma}
\begin{pf}
Take any $u\in\dom^{\bar D}\cap C_0(\overline D)$ such that
$\operatorname{supp}[u]\subset G$. From (\ref{equ6.3}) and (\ref{equ6.4}), we then have
%
%
\begin{equation}\label{equ6.13}
\eta_D(u,v)=-\int_{G\times(\bar D\setminus\bar
G)}u(y)v(x)k(x,y)m_D(dx)m_D(dy).
\end{equation}
We can now proceed as in \cite{FOT94}, page 163. The function $g_v$
defined by
(\ref{equ6.12}) belongs to $L^2(G;m_G)$ on account of   condition
(\ref{cond-1})
on the kernel $k$. Therefore, we obtain from
(\ref{equ6.13})
\begin{eqnarray*}
\eta_{G,\alpha}^0(R_\alpha^{G,0}g_v, u)&=&\int_G g_v(x) u(x) m_G(dx)\\
&=&\int_{G\times(\bar D\setminus\bar G)}u(x)v(y)k(x,y)m_D(dx)m_D(dy)\\
&=&-\eta_D(v,u)=-\eta_{D,\alpha}(v,u)\\
&=&-\eta_{G,\alpha}^0(v-H_\alpha^{\bar D\setminus G}v, u),\qquad \alpha
>\beta_0,
\end{eqnarray*}
the last identity being a consequence of (\ref{equ6.10}). Since $\dom
^{\bar
D}\cap C_0(G)$
is $\eta_{G,\alpha}^0$-dense in $\dom_G^0$, we get
\[
H_\alpha^{\bar D\setminus G}v(x) =H_\alpha^{\bar D\setminus G}v(x)
-v(x)=R_\alpha^{G,0}g_v(x) \qquad\mbox{for } m_G\mbox{-a.e. on } G.
\]
We then obtain (\ref{equ6.11}) because $H_\alpha^{\bar D\setminus G}v$ and
$R_\alpha^{G,0}g_v$ are $\eta^0_G$-quasi-continuous by (\ref{equ6.10}).
\end{pf}
\begin{theorem}\label{theo6.2}
\begin{longlist}[(iii)]
\item[(i)]
$X^{\bar D}=(X_t, P_x)$ admits no jump from
$D$ to
$\partial D$:
%
%
\begin{equation}\label{equ6.14}
P_x(X_{t-}\in D, X_t\in\partial D \mbox{ for some } t>0)=0\qquad
\mbox{for q.e. } x\in D.\hspace*{-25pt}
\end{equation}

\item[(ii)] If $D$ is relatively compact, then $X^{\bar D}$ is conservative:
denoting by $\zeta$ the lifetime of $X^{\bar D}$,
%
%
\begin{equation}\label{equ6.15}
P_x(\zeta=\infty)=1 \qquad\mbox{for q.e. } x\in\overline D.
\end{equation}

\item[(iii)] If $D$ is relatively compact, then $X^{D,0}=(X_t^{D,0},
P_x)$ admits no
killing inside~$D$: denoting by $\zeta^0$ the lifetime of $X^{D,0}$,
%
%
\begin{equation}\label{equ6.16}
P_x(X^{D,0}_{\zeta^0-}\in D, \zeta^0<\infty)=0 \qquad\mbox{for
q.e. } x\in D.
\end{equation}
\end{longlist}
\end{theorem}
\begin{pf}
(i) For any open set $G$ as Lemma \ref{lemma6.1} and any compact
subset $F$ of
$\partial D$, we can find a uniformly bounded sequence
$\{v_n\}\subset\dom^{\bar D}\cap C_0(\overline D)$ with support being
contained in
a common compact subset of $\overline D\setminus\overline G$ and
$\lim_{n\to\infty}v_n=1_F$. Then $g_{v_n}(x)$ are uniformly bounded
and converge to
$g_{1_F}(x)=0$ as $n\to\infty$. Therefore, by letting $n\to\infty$
in (\ref{equ6.11}) with
$v_n$ in place of $v$, we get $P_x(X_{\tau_G}\in F
)=0$ for q.e.
$x\in G$. Since $G$ and $F$ are arbitrary with the stated properties,
we have (\ref{equ6.14}).

\mbox{}\phantom{i}(ii) When $D$ is relatively compact, $1\in C_0^{\mathrm{lip}}(\overline D)$
so that we
see from (\ref{equ6.3}) and (\ref{equ6.4}) that $1\in\dom^{\bar D}$ and
$\eta
_D(1,v)=0$ for any
$v\in\dom^{\bar D}$. We have therefore, for any $\alpha>\beta_0$ and
$f\in L^2(\overline D,m_D)$,
\[
0=\eta_D(1,\widehat G_\alpha f)=(1, f)_D-\alpha(1, \widehat G_\alpha f)_D
=(1-\alpha R_\alpha1, f)_D,
\]
where $\widehat G_\alpha$ is the dual resolvent. This implies that
$\alpha R_\alpha1=1$ $m_D$-a.e. for $\alpha>\beta_0$ and
consequently q.e. on
$\overline D$ because $R_\alpha1$ is quasi-continuous. Equation~(\ref
{equ6.15}) is proven.

(iii) This is an immediate consequence of (i), (ii) as $X^{D,0}$ is
the part process of $X^{\bar D}$ on $D$.
\end{pf}

We conjecture that the property (\ref{equ6.16}) for $X^{D,0}$ holds
true without the
assumption of the relative compactness of $D$ and especially for the minimal
process $X^0$ on $E$.

Finally, we consider the case where $E$ is $\rd$ and $m$ is the
Lebesgue measure
on it. For $\alpha\in(0,2)$ and an arbitrary open set $D\subset\rd
$, we make
use of the L\'{e}vy kernel
\[
k^{[\alpha]}(x,y)=
\frac{\alpha2^{\alpha-1}\Gamma(({\alpha+d})/{2})}{\pi^{d/2}
\Gamma(1-\alpha/2)} \frac{1}{|x-y|^{d+\alpha}},\qquad x,y \in\rd,
\]
of the symmetric $\alpha$-stable process to introduce the Dirichlet form
%
%
\begin{equation}\label{equ6.17}
\cases{
\displaystyle \form_D^{[\alpha]}(u,v):= \iint_{D\times D\setminus\di}
\bigl(u(y)-u(x)\bigr)
\bigl(v(y)-v(x)\bigr)k^{[\alpha]}(x,y) \,dx\,dy,\vspace*{2pt}\cr
\displaystyle \dom^{[\alpha],r}_D=\bigl\{u\in L^2(D)\dvtx u \mbox{ is Borel
measurable and }
\form_D^{[\alpha]}(u,u)<\infty\bigr\},}\hspace*{-40pt}
\end{equation}
on $L^2(D)$ based on the Lebesgue measure on $D$. Denote by
$\dom_{\bar D}^{[\alpha]}$ the $\form^{[\alpha]}_{D, 1}$-closure\vadjust{\goodbreak} of
$C_0^{\mathrm{lip}}(\overline D)$ in $\dom_D^{[\alpha],r}$. For $s\in(0,d]$,
a Borel subset $\Gamma$ of $\rd$ is said to be an $s$-\textit{set} if
there exist
positive constants $c_1, c_2$ such that for all $x\in\Gamma$ and
$r\in(0,1]$,
$c_1 r^s\le\mathcal{H}^s(\Gamma\cap B(x,r))\le c_2r^s$, where
$\mathcal{H}^s$ denotes
the $s$-dimensional Hausdorff measure on $\rd$ and $B(x,r)$ is the
ball of radius
$r$ centered at $x\in\rd$.

If the open set $D$ is a $d$-set, then, by making use of
Jonsson--Wallin's trace
theorem \cite{JW84} as in \cite{BBC03}, one can show that
$\dom_{\bar D}^{[\alpha]}=\dom_D^{[\alpha],r}$ and moreover that a~subset of
$\overline D$ is $\form_D^{[\alpha]}$-polar iff it is polar with
respect to the
symmetric $\alpha$-stable process on~$\rd$.

Let us consider the kernel $k^{(1)}$ of (\ref{stable-like-kernel}) for
$w(x)$ given by (\ref{weight}) and
$\alpha(x)$ satisfying condition (\ref{stable-like}). In particular, it
is assumed that
\[
0<\underline\alpha\le\alpha(x)\le\overline\alpha<2
\]
for some constant $\underline\alpha, \overline\alpha$. $k^{(1)}$
satisfies
conditions (\ref{cond-1})--(\ref{cond-4}) by Proposition~\ref{prop5.1}
and one can associate with it
the regular
lower bounded semi-Dirichlet form~$\eta_D$ (resp., $\eta_D^0$) on
$L^2(\overline D; 1_D \,dx)$ [resp., $L^2(D)$] possessing as its
reference form~$\form_D$ (resp.,~$\form_D^0$) defined right after (\ref{equ6.1}) for $k^{(1)}$
and the
Lebesgue measure in place of $k$ and~$m$.

Suppose $D$ is bounded, then there exist positive constants $c_3, c_4$ with
\[
c_3 k^{[\underline\alpha]}(x,y)\le k_s^{(1)}(x,y) \le
c_4 k^{[\overline\alpha]}(x,y),\qquad x,y\in\overline D,
\]
so that
%
%
\begin{equation}\label{equ6.18}
c_3 \form^{[\underline\alpha]}_D(u,u)\le\form_D(u,u)\le c_4
\form^{[\overline\alpha]}_D(u,u),\qquad u\in C_0^{\mathrm{lip}}(\overline D).
\end{equation}

For the kernel $k^{(1)}$, the Hunt process $X^{\bar D}$ on $\overline
D$ associated
with $(\eta_D, \dom_{\bar D})$\vspace*{2pt} is called a modified reflecting
stable-like process,
while its part process $X^{D,0}$ on $D$, which is associated with
$(\eta_D^0, \dom_D^0)$, is called a censored stable-like
process. 
\begin{prop}\label{prop6.1}
Assume that $D$ is a bounded open $d$-set.

\begin{longlist}[(ii)]
\item[(i)]
If $\partial D$ is polar with respect to the symmetric
$\overline\alpha$-stable process on $\rd$, then the censored
stable-like process
$X^{D,0}=(X_t^{D,0}, P_x, \zeta^0)$ is conservative and it does not
approach to
$\partial D$ in finite time:
%
%
\begin{equation}\label{equ6.19}
P_x(\zeta^0=\infty)=1,\qquad P_x(X_{t-}^{D,0}\in\partial D
\mbox{ for some } t>0)=0.
\end{equation}

\item[(ii)] If $\partial D$ is nonpolar with respect to the symmetric
$\underline\alpha$-stable process on~$\rd$, then the censored stable-like
process $X^{D,0}$ satisfies
%
%
\begin{equation}\label{equ6.20}
\int_D P_x(X_{\zeta^0-}^{D,0}\in\partial D, \zeta^0<\infty)h(x)\,dx=
\int_D P_x(\zeta^0<\infty)h(x)\,dx>0\hspace*{-28pt}
\end{equation}
for any strictly positive Borel function $h$ on $D$ with $\int_Dh(x)\,dx=1$.
\end{longlist}
\end{prop}
\begin{pf}
(i) Since $\form_D$ is a reference form of $(\eta_D,
\dom_{\bar D})$,
we see that $\partial D$ is $\eta_D$-polar by (\ref{equ6.18}) and the stated
observation in
\cite{BBC03}. The assertions of (i) then follows from Theorem \ref
{theo4.2}(ii) and
Theorem 6(ii).\vadjust{\goodbreak}

(ii) $\partial D$ is not $\eta_D$-polar by (\ref{equ6.18}) and
accordingly not
$m$-polar
with respect to the process $X^{\bar D}$ by Theorem \ref{theo4.2}(v),
where $m$
is the Lebesgue
measure on~$D$. Taking Theorem \ref{theo6.2}(i), (iii) into account, we
then get (\ref{equ6.20}).
\end{pf}

The polarity of a set $N\subset\rd$ with respect to the symmetric
$\alpha$-stable process
is equivalent to $C^{\alpha/2,2}(N)=0$ for the Bessel capacity
$C^{\alpha/2,2}$ (cf. Section 2.4
of the second edition of \cite{FOT94}). The latter has been well
studied in \cite{AH96}
in relation to the Hausdorff measure and the Hausdorff content.
For instance, when $\alpha\le d$ and $\partial D$ is a $s$-set,
$\partial D$ is polar in this sense if and only if $\alpha+s\le d$.
Of course, we get the same results as above for the second kernel
$k^{(1)*}$ in (\ref{stable-like-kernel}).

%

%
\printaddresses

\end{document}